\documentclass{article}
\usepackage{fullpage}
\usepackage{amsmath,amsthm,mathrsfs, verbatim,color,lscape,xspace}
\usepackage{listings}
\usepackage{graphicx}
\usepackage{secdot}
\usepackage{amsmath, dsfont}
\usepackage{amsfonts}
\usepackage{amssymb}
\usepackage{enumerate}
\usepackage{ifpdf}
\usepackage{lscape}
\usepackage[tight]{subfigure}

\newcommand{\erfct}[1]{\frac{2}{\sqrt\pi}\int_{#1}^\infty e^{-x^2}dx}

\newcommand{\lt}[1]{\ensuremath{\hat{#1}}}          

\newcommand{\apr}[1]{\ensuremath{\tilde{#1}}}
\newcommand{\com}[1]{\ensuremath{\bar{#1}}}

\newcommand{\disf}[2][u]{\ensuremath{{#2}\left(#1\right)}}     
\newcommand{\df}[2][u]{\ensuremath{{#2}(#1)}}       
\newcommand{\distance}[1]{\ensuremath{\mathcal{D}(#1)}}
\newcommand{\supnorm}[3][t]{\ensuremath{\textup{sup}_{\substack{#1}}\left|#2(#1) -#3(#1)\right|}}

\newcommand{\con}[3][n]{\ensuremath{{#2}^{*#1}_{#3}}}            
\newcommand{\cona}[3][n]{\ensuremath{\tilde{{#2}}^{*#1}_{#3}}}   

\newcommand{\e}[1][]{\ensuremath{\mathbb{E}{#1}}\xspace}    
\newcommand{\pr}[1][]{\ensuremath{\mathbb{P}{#1}}\xspace}   
\newcommand{\ind}[3]{\ensuremath{{#1}^{#2}_{#3}}}           

\newcommand{\footnoteremember}[2]
{
   \newcounter{#1}\footnote{#2}\setcounter{#1}{\value{footnote}}
}
\newcommand{\footnoterecall}[1]
{
   \footnotemark[\value{#1}]
}

\newtheorem{proposition}{Proposition}
\newtheorem{lemma}{Lemma}
\newtheorem{definition}{Definition}
\newtheorem{theorem}{Theorem}
\theoremstyle{definition}
\newtheorem{remark}{Remark}

\begin{document}

\title{On the accuracy of phase-type approximations\\ of heavy-tailed risk models}
\author{
     E. Vatamidou\footnoteremember{TU/e}{Depart. of Mathematics \& Computer Science, Eindhoven University of Technology, P.O. Box 513, 5600 MB Eindhoven, The Netherlands}\footnoteremember{EURANDOM}{\textsc{Eurandom}, Eindhoven University of Technology, P.O. Box 513, 5600 MB Eindhoven, The Netherlands}\\
     \small \texttt{e.vatamidou@tue.nl}\\
     \and
     I.J.B.F. Adan\footnoterecall{EURANDOM}\footnoteremember{MechEng}{Dept.\ of Mechanical Engineering, Eindhoven University of Technology, P.O. Box 513, 5600 MB Eindhoven, The Netherlands}\\
     \small \texttt{i.j.b.f.adan@tue.nl}\\
     \and
     M. Vlasiou\footnoterecall{TU/e}\footnoterecall{EURANDOM}\\
     \small \texttt{m.vlasiou@tue.nl}
     \and
     A.P. Zwart\footnoterecall{EURANDOM}\footnoteremember{CWI}{Centrum Wiskunde \& Informatica (CWI), P.O. Box 94079, 1090 GB Amsterdam, The Netherlands}\\
     \small \texttt{Bert.Zwart@cwi.nl}\\
}
\maketitle

\begin{abstract}
  Numerical evaluation of ruin probabilities in the classical risk model is an important problem. If claim sizes are heavy-tailed, then such evaluations are challenging. To overcome this, an attractive way is to approximate the claim sizes with a phase-type distribution. What is not clear though is how many phases are enough in order to achieve a specific accuracy in the approximation of the ruin probability. The goals of this paper are to investigate the number of phases required so that we can achieve a pre-specified accuracy for the ruin probability and to provide error bounds. Also, in the special case of a completely monotone claim size distribution we develop an algorithm to estimate the ruin probability by approximating the excess claim size distribution with a hyperexponential one. Finally, we compare our approximation with the heavy traffic and heavy tail approximations.
\end{abstract}

\section{Introduction}\label{Intro}

In this paper we deal with the numerical evaluation of the ruin probability in the classical risk model \cite{prabhu61}. In this model, we have claims (for money) which arrive to an insurance company according to a Poisson process. The probability $\df{\psi}$ of ultimate ruin is the probability that the reserve of the insurance company ever drops below zero, where $u$ is the initial capital of the company and where the total income (premium) rate is 1.

In many financial applications, the classical assumption of exponentially decaying claim sizes is not applicable \cite{embrechts-MEE}. An appropriate way to model claim sizes in such cases is by using heavy-tailed distributions. Such distributions decay more slowly than any exponential function, which means that with such distributions there exists a nontrivial probability of an extremely large claim size \cite{asmussen-APQ,rolski-SPIF}.

Heavy-tailed distributions also play a significant role in queueing models, where service times can take extremely large values. It is actually well known \cite{asmussen-RP,asmussen-APQ} that the probability of eventual ruin for an insurance company with an initial cash reserve $u$ is equal to the stationary waiting probability $\pr (W_q>u)$ of a G/G/1 queue, where service times in the queueing model correspond to the random claim sizes.

In this paper, we assume that claim sizes arrive according to a Poisson process. Therefore, the ruin probability can be found by using the well-known Pollaczek-Khinchine formula \cite{asmussen-RP}. This formula involves the convolutions of the excess claim size distribution (see Section~\ref{error bound}), which cannot be easily computed, and thus one usually resorts to Laplace transforms. However, a major difficulty when analyzing models with heavy-tailed distributions is that Laplace transforms of such distributions oftentimes do not have an analytic closed form. This is, in particular, the case for the Pareto and Weibull distributions. Thus, analytic methods, which use the Laplace transform of the claim sizes, are difficult \cite{abate99b} or even impossible to use in such cases.

When the ruin probability, or equivalently the waiting time distribution of a G/G/1 queue, cannot be computed exactly it needs to be approximated. The approximations for the ruin probability can be classified in two general categories: based on the average amount of claims per unit time (or load of the system) and based on the characteristics of the claim size distribution.

In the first category we have the heavy traffic \cite{kalashnikov-GSBREA,kingman62} and light traffic \cite{asmussen92b,bloomfield72,daley84,daley91,sigman92} approximations. If, on average, the premiums exceed only slightly the expected claims then most appropriate for modeling is the heavy traffic approximation. The drawback of this approximation though is that it requires finite first two moments for the claim size distribution, a condition which may not be satisfied for several heavy-tailed distributions. On the other hand, when on average, the premiums are much larger than the expected claims then the light traffic approximation is used. In risk theory, heavy traffic is most often argued to be the typical case rather than light traffic, which makes the light traffic approximation only of limited interest.

Closely related to the previous approximations is the Edgeworth series expansion \cite{wallace58}, which is a refinement of the central limit theorem. Asymptotic results for the ruin probability are given in \cite{blanchet10} and these approximations can be useful in applications where moments are computable, but the distribution is not.

As mentioned above, another category of approximations is based on the characteristics of the claim size distribution. Two known approximations, which are based on the idea of matching the moments of the ruin probability, are the Beekman-Bower's \cite{beekman69} and the de Vylder's \cite{vylder78} approximations. However, for some heavy-tailed distributions, such as the Pareto mentioned above, (higher-order) moments may be infinite, thus making conventional moment-matching methods fail.

A particularly effective approach in handling distributions with infinite moments is the Transform Approximation Method (TAM). The Laplace transform of a positive definite distribution, like the claim size distribution, exists always even if it does not have a closed analytic form. The TAM is based on the idea of approximating the Laplace transform of the claim sizes rather than directly their distribution \cite{harris00,harris98,shortle04}. A drawback of this method though is that the accuracy of the approximation of the ruin probability cannot be predetermined.

When the claim sizes belong in the subexponential class of distributions \cite{teugels75} then the heavy tail \cite{bahr75,borovkov92,embrechts82,pakes75} approximation is also used. However, the disadvantage of this approximation is that it provides a good fit only at the tail of the ruin probability, especially in the case where the average amount of claim per unit time is close to one.

Finally, a natural approach to address the problem of non-existence of the Laplace transform for a heavy-tailed distribution in a closed form is to approximate the claim size distribution with a phase-type distribution \cite{feldmann98,lucantoni94,neuts-SSM}. The main advantage of approximating a heavy-tailed claim size distribution with a phase-type distribution is that in the latter case the Laplace transform of the claim sizes has a closed form. Several approximation methods for probability distributions using special cases of phase-type distributions, such as the Coxian and the hyperexponential distributions, have been proposed \cite{asmussen96,sasaki04,starobinski00}. These methods can provide approximations for the claim sizes with high accuracy. However, one of their disadvantages is that the accuracy of the approximation of the ruin probability cannot be predetermined. Another drawback is that the number of phases needed to achieve a desired accuracy cannot be determined a priori and most times it is found by trial and error.

In this paper, we develop a new approach for approximating the ruin probability, when the claim sizes follow a heavy-tailed distribution. From the Pollaczek-Khinchine formula (see Section~\ref{error bound}) it is clear that in order to evaluate the ruin probability, we only need to have a closed analytic form for the Laplace transform of the excess claim size distribution. For this reason, instead of approximating the claim size distribution, we approximate directly the stationary excess distribution with a hyperexponential distribution, a special case of a phase-type distribution. Since the Laplace transform of a hyperexponential distribution exists in a closed analytic form, we can numerically evaluate $\df{\psi}$ by inverting its Laplace transform.

An advantage of our approximation, which we call the {\it spectral approximation}, is that it has a predetermined accuracy. Thus, we first choose the accuracy we want to achieve in our approximation, and later on we determine the number of states for the hyperexponential distribution that are sufficient to guarantee this accuracy. Another interesting feature is that the bound that we guarantee is valid for the whole domain of the ruin probability $\df{\psi}$ and not only for a subset of it, contrary to other bounds that exist in the literature \cite{kalasnikov99,starobinski00}.

In Section~\ref{error bound}, we find bounds for the $n$th convolution of the excess claim size distribution. We prove that the bound for the convolution is linear with respect to the accuracy we choose for the excess claim size distribution. We also give the main result of this paper, which is the error bound for the ruin probability $\df{\psi}$. Later, we focus on a class of heavy-tailed distributions that are in addition completely monotone, and we show that we can always approximate a completely monotone distribution with a hyperexponential one for any desired accuracy. We also prove that if the claim size distribution is completely monotone with finite mean, then the stationary excess distribution is also completely monotone. Finally, we sketch {\it the spectral approximation algorithm}, which approximates a completely monotone excess claim size distribution with a hyperexponential distribution for any desired accuracy.

Later on, we also compare the spectral approximation with the heavy traffic and the heavy tail approximations. Thus, in Section~\ref{heavy-traffic and heavy-tail approximations}, we give the basic characteristics of the latter two approximations, and mention their advantages and disadvantages. We devote Section~\ref{numerical results} to numerical results. We do a series of experiments in order to compare the spectral approximation with the heavy traffic and the heavy tail approximations. As test distributions we use the Pareto, the Weibull and a class of long-tail distributions introduced in \cite{abate99a}. In addition, we extend a bound that is given in the literature \cite{brown90} for the heavy traffic approximation to a specific case of the heavy traffic approximation that we use in our experiments.

In Section~\ref{conclusions}, we discuss the results. Finally, in the Appendix we present the steps of the algorithm that we use to produce the spectral approximation.

\section{Spectral approximation for the ruin probability}\label{error bound}
Consider the classical compound Poisson risk model \cite{asmussen-RP}. In broad terms, a risk reserve process $\{R_t\}_{t \geq 0}$ is a model for the time evolution of the reserve of an insurance company, where the initial reserve is denoted by $u=R_0$. Claims arrive according to a Poisson process $\{N_t\}_{t \geq 0}$ with rate $\lambda$. The claim sizes $U_1,U_2,\dots$ are i.i.d. with common distribution $B$ and independent of $\{N_t\}$, and premiums flow in at a rate 1 per unit time. Putting all these together we see that
\begin{equation*}\label{risk reserve process}
  R_t = u+t-\sum_{k=1}^{N_t}U_k.
\end{equation*}
For mathematical purposes, it is frequently more convenient to work with the claim surplus process $\{S_t\}_{t \geq 0}$ which is defined as $S_t =u - R_t$; as one can see from the expression above, this is merely a compound Poisson process with positive jumps and negative drift, a process well studied in the literature. The probability $\df{\psi}$ of ultimate ruin is the probability that the reserve ever drops below zero, or equivalently the probability that the maximum $M = \sup_{0 \leq t < \infty} S_t$ ever exceeds $u$; i.e.
\begin{equation}
  \df{\psi} = \pr(M>u).
\end{equation}

Since we consider Poisson arrivals for the claims, for the evaluation of the ruin probability, the well-known Pollaczek-Khinchine formula \cite{asmussen-RP} can be used:
\begin{equation}\label{polllaczek khinchine formula}
  \df{\psi} = (1-\rho)\sum_{n=0}^\infty\rho^n \disf{\overline{{\con{B}{0}}}},
\end{equation}
where $\rho <1$ is the average amount of claim per unit time. For a distribution $F$ we use the notation $\con{F}{}$ to denote its $n$th convolution, $\com{F}$ to denote its complementary cumulative distribution (i.e. the tail) and $\lt F$ to denote its Laplace transform. Moreover, $B_0$ is the stationary excess distribution, which is defined as
\begin{equation*}
  \disf{B_0} = \frac{1}{\e U} \int_0^u \disf[t]{\com{B}}dt,
\end{equation*}
where $\e U$ is the (finite) mean of the claim sizes. The $n$th moment of the claim sizes is denoted by $\e U^n$.

For the evaluation of $\df{\psi}$, \eqref{polllaczek khinchine formula} is not entirely satisfying because the infinite sum of convolutions at the right-hand side of the formula cannot be easily computed analytically and sometimes not even numerically. In order to overcome this difficulty we use Laplace transforms, which convert convolutions of distributions into powers of their Laplace transform. In terms of Laplace transforms, the Pollaczek-Khinchine formula can be written as:
{
\begin{equation}\label{laplace transform of polllaczek khinchine formula}
  \e e^{-s M} = (1-\rho)\sum_{n=0}^\infty\rho^n \disf[s]{\lt{B}_0^n} = \frac{1-\rho}{1 - \rho \disf[s]{\lt B_0}}.
\end{equation}
From \eqref{laplace transform of polllaczek khinchine formula} it is clear why it is necessary to have a closed analytic form only for the Laplace transform of the {\it excess} claim size distribution, rather than the claim size distribution itself. Thus, the main idea of our algorithm is to approximate the excess claim size distribution with a phase-type distribution, which has a closed analytic Laplace transform, and apply Laplace inversion to evaluate the ruin probability.

\subsection{Error bound for the ruin probability}

In this section, we provide a bound for the ruin probability when we approximate the excess claim size distribution with a known distribution, e.g. a phase-type distribution. If we approximate $B_0$ with a known distribution (not only a phase-type) then we can compute the ruin probability through the Pollaczek-Khinchine formula \eqref{polllaczek khinchine formula}. From \eqref{polllaczek khinchine formula} and the triangular inequality, the error between the ruin probability and its approximation is then
\begin{equation}\label{ruinerror}
    \left|\df{\psi} - \df{\apr{\psi}} \right| = \left|\sum_{n=0}^\infty(1-\rho)\rho^n  \left(\disf{\con{B}{0}}-\disf{\cona{B}{0}}\right)\right|
                                            \le \sum_{n=0}^\infty(1-\rho)\rho^n \left| \disf{\con{B}{0}}-\disf{\cona{B}{0}}\right|,
\end{equation}
where $\apr{F}$ denotes the approximation of a distribution $F$, and $\apr{\psi}$ is the exact result we obtain from the Pollaczek-Khinchine formula for the ruin probqbility when we use an approximate claim size distribution. From \eqref{ruinerror} we see that as a first step to find a bound for the ruin probability is to find a bound for the difference $\left| \disf{\con{B}{0}}-\disf{\cona{B}{0}} \right|$. This is given in the following proposition.

\begin{proposition}\label{stationaryexcessbound}
  If $\supnorm[x]{\ind{B}{}{0}}{\ind{\apr{B}}{}{0}}\le \epsilon$ for $x\in [0,u]$, then $\left| \disf{\con{B}{0}}-\disf{\cona{B}{0}}\right|\le~ n\epsilon$.
\end{proposition}
\begin{proof}
  We prove this by induction. For $n=2$,
\begin{align*}
   |\disf[u]{\con[2]{B}{0}}-\disf[u]{\cona[2]{B}{0}}| &= |\disf{\ind{B}{}{0}*\ind{B}{}{0}} \pm
                                              \disf{\ind{\apr{B}}{}{0}*\ind{B}{}{0}} -
                                              \disf{\ind{\apr{B}}{}{0}*\ind{\apr{B}}{}{0}}|\\
                                          &\le |\disf{{(\ind{B}{}{0}-\ind{\apr{B}}{}{0})}*\ind{B}{}{0}}| +
                                             |\disf{(\ind{B}{}{0}-\ind{\apr{B}}{}{0})*\ind{\apr{B}}{}{0}}|\\
                                          &\le \int_0^u\underbrace{|\disf[u-x]{(\ind{B}{}{0}-\ind{\apr{B}}{}{0})}|}_{\le\epsilon}d\disf[x]{\ind{B}{}{0}} +
                                          \int_0^u\underbrace{|\disf[u-x]{(\ind{B}{}{0}-\ind{\apr{B}}{}{0})}|}_{\le\epsilon}d\disf[x]{\ind{\apr{B}}{}{0}}\\
                                          &\le \epsilon \disf{\ind{B}{}{0}} + \epsilon \disf{\ind{\apr{B}}{}{0}}\\
                                          &\leq 2\epsilon.
   \intertext{Assume now that the bound holds for a fixed $n$. We prove that it also holds for $n+1$.}
|\disf{\con[(n+1)]{B}{0}}-\disf{\cona[(n+1)]{B}{0}}|  &= |\disf{\ind{B}{}{0}*\con{B}{0}} \pm
                                                            \disf{\ind{\apr{B}}{}{0}*\con{B}{0}}
                                                            -\disf{\ind{\apr{B}}{}{0}*\cona{B}{0}}|\\
                                                      &\le |\disf{(\ind{B}{}{0}-\ind{\apr{B}}{}{0})*\con{B}{0}}|
                                                            +|\disf{\ind{\apr{B}}{}{0}*(\con{B}{0}-\cona{B}{0})}|\\
                                                      &\le \int_0^u\underbrace{|\disf[u-x]{(\ind{B}{}{0}-\ind{\apr{B}}{}{0})}|}
                                                            _{\le\epsilon} d\disf[x]{\con{B}{0}}
                                                            +\int_0^u\underbrace{|\disf[u-x]{(\con{B}{0}-\cona{B}{0})}|}
                                                            _{\le n\epsilon} d\disf[x]{\ind{\apr{B}}{}{0}}\\
                                                      &\le\epsilon \disf{\con{B}{0}} + n\epsilon \disf{\cona{B}{0}}\\
                                                      &\leq (n+1)\epsilon.
\end{align*}
\end{proof}

In words, Proposition~\ref{stationaryexcessbound} says that if we bound the excess claim size distribution with some accuracy $\epsilon$, then a bound for its $n$th convolution is linear with respect to this accuracy $\epsilon$. Consequently, from Proposition~\ref{stationaryexcessbound}, we have the following result.

\begin{proposition}\label{ruinbound}
  If $\supnorm[x]{\ind{B}{}{0}}{\ind{\apr{B}}{}{0}}\le \epsilon$ for $x\in [0,u]$, then a bound for the ruin probability is
  \begin{equation*}\label{ruinbound}
    \left|\df{\psi}-\disf{\apr{\psi}} \right|\leq \frac{\epsilon\rho}{1-\rho}.
  \end{equation*}
\end{proposition}
\begin{proof}
  \begin{align*}
    \left|\df{\psi}-\df{\apr{\psi}}\right|  &\le \sum_{n=0}^\infty(1-\rho)\rho^n \left|\disf{\con{B}{0}}-\disf{\cona{B}{0}}\right|\\
                                            & \le \sum_{n=0}^\infty(1-\rho)\rho^n n\epsilon
                                             =\epsilon\rho(1-\rho)\sum_{n=0}^\infty n\rho^{n-1}\\
                                            &=\epsilon\rho(1-\rho)\left(\frac{1}{1-\rho}\right)^\prime
                                             =\epsilon\rho(1-\rho)\frac{1}{(1-\rho)^2}\\
                                            &=\frac{\epsilon\rho}{1-\rho}.
  \end{align*}
\end{proof}

Notice that the bound in Proposition~\ref{ruinbound} is independent of $u$, for all $u\geq 0$. Thus, if we define the sup norm distance between two positive definite distributions $F_1$ and $F_2$ as $\mathcal{D}(F_1,F_2)=\textup{sup}_u|\disf{F_1}-\disf{F_2}|$, $u\geq0$, we conclude that $\distance{\psi,\apr{\psi}} \leq \frac{\epsilon\rho}{1-\rho}$, whenever $\distance{B_0,\apr{B}_0} \leq \epsilon$. Observe that the term $1-\rho$ at the denominator has as consequence that, higher load $\rho$ requires a more accurate approximation of the $B_0$ to obtain tight bounds for the ruin probability.

To sum up, when the excess claim size distribution is approximated with some desired accuracy $\epsilon$, then a bound for the ruin probability, which is linear with respect to $\epsilon$, is guaranteed by Proposition~\ref{ruinbound}. Thus, our next goal is to develop a way to approximate the excess claim size distribution with a hyperexponential one, a particular case of a phase-type distribution, with any desired accuracy. We complete this step in the next section.

\subsection{Completely monotone claim sizes}\label{completely monotone claim sizes}
We are interested in evaluating the ruin probability when the claim sizes follow a heavy-tailed distribution, such as Pareto or Weibull. These two distributions belong also to the class of completely monotone distributions, which is defined below.
\begin{definition}
  A probability density function (pdf) is said to be completely monotone (c.m.) if all derivatives of f exist and if
  \begin{equation*}\label{definition of complete monotonicity}
    (-1)^n \df{\ind{f}{(n)}{}} \geq 0 \text { for all } u>0 \text{ and } n\geq 1.
  \end{equation*}
\end{definition}

Completely monotone distributions can be approximated arbitrarily close by hyperexponentials \cite{feldmann98}. Here, we provide a method to approximate a completely monotone excess claim size distribution with a hyperexponential one in order to achieve any desired accuracy for the ruin probability. The following result is standard; see e.g. \cite{feller-IPTIA}.

\begin{theorem}\label{formula for completely monotone distributions}
  A pdf is called completely monotone if and only if it is a mixture of exponential pdf's. That is,
  \begin{equation*}
    \df{f} =\int_0^{+\infty}y e^{-y u}d\disf[y]{G}, \ u\geq0,
  \end{equation*}
  for some proper positive-definite cumulative distribution function (cdf) G. We call G the \textup{spectral} cdf. For the tail or the complementary cumulative distribution function (ccdf) of a completely monotone distribution it holds that
  \begin{equation*}
    \disf{\com{F}}   =\int_u^{+\infty}\df[x]{f}dx
                         =\int_0^{+\infty}\int_u^{+\infty}y e^{-y x}dxd\disf[y]{G}
                         =\int_0^{+\infty}e^{-y u}d\disf[y]{G}.
  \end{equation*}
\end{theorem}

An alternative way to define a c.m. distribution is by using Laplace transforms. From Theorem~\ref{formula for completely monotone distributions} it is obvious that a pdf $f$ is c.m. if its tail can be written as the Laplace transform of some positive-definite distribution $G$. The following lemma is an immediate consequence.

\begin{lemma}\label{the complementaty integrated tail}
  If the claim size distribution is c.m. then the excess claim size distribution is c.m. too.
\end{lemma}
\begin{proof}
  If $B$ is a completely monotone distribution, then $\disf{\com{B}}=\int_0^{+\infty}e^{-yu}d\disf[y]{G}$, for some spectral function $G$. Thus,
  \begin{align*}
    \disf{\com{B}_0}  &=\frac{1}{\e U}\int_u^{+\infty}\disf[x]{\com{B}}dx
                       =\frac{1}{\e U}\int_u^{+\infty}\int_0^{+\infty}e^{-y x}d\disf[y]{G}dx
                       =\frac{1}{\e U}\int_0^{+\infty}d\disf[y]{G} \int_u^{+\infty}e^{-y x} dx\\
                      &=\int_0^{+\infty}e^{-y u}\frac{d\disf[y]{G}}{y\e{U}}
                       =\int_0^{+\infty}e^{-y u}d\disf[y]{H},
  \end{align*}
  where $d\disf[y]{H}=\frac{d\disf[y]{G}}{y\e{U}}$.
\end{proof}

In this paper, we are interested in finding a bound for the excess claim size distribution. In order to achieve our goal, we approximate the spectral function of the excess claim size distribution by a step function with some fixed (and pre-determined) accuracy $\epsilon$ and then calculate the error of the approximation for the excess claim size distribution itself.

\begin{lemma}
  Let $G$ be the spectral function  of the c.m. excess claim size distribution $B_0$, and let the step function $\apr{G}$ satisfy $\distance{G,\apr{G}} \leq \epsilon$. Then, $\distance{B_0,\apr{B}_0} \leq \epsilon$, where $\apr{B}_0$ is the c.m. distribution with spectral function $\apr{G}$.
\end{lemma}
\begin{proof}
  Since the spectral cdf $G$ is proper, we have by definition that it has no atom at 0 and that it is right continuous. Thus, $\disf[0]{G}=0$ and $\disf[+\infty]{G}=1<\infty$. Then it holds that
  \begin{align*}
    \int_0^{+\infty}e^{-u y}d\disf[y]{G} &=e^{-u y}\disf[y]{G}\mid^{+\infty}_0-\int_0^{+\infty}\disf[y]{G}de^{-u y}\\
                                         &=\int_0^{+\infty}u e^{-u y}\disf[y]{G}dy.
  \end{align*}
Suppose now that $ \distance{G,\apr{G}} \leq \epsilon$. Then
\begin{align*}
  \left| \disf{\com{B}_0}-\disf{\apr{\com{B}}_0}\right|   &\leq \left| \int_0^{+\infty}(\disf[y]{G}-
                                                                    \disf[y]{\apr{G}})u e^{-u y}dy \right|\\
                                                          &\leq \int_0^{+\infty}\underbrace{\left|\disf[y]{G}
                                                                    -\disf[y]{\apr{G}}\right|}_{\leq \epsilon}u e^{-u y}dy \leq \epsilon,
\end{align*}
for all $u\geq0$. So, $\distance{B_0,\apr{B}_0} \leq \epsilon$.
\end{proof}

Summarizing, if we want to approximate the claim size distribution with a hyperexponential\footnote{By definition, a hyperexponential distribution with $k$ phases is a c.m. distribution with spectral function a step function with $k$ jumps.} with some fixed accuracy $\epsilon$, it is sufficient to approximate the spectral cdf of the c.m.\ excess claim size distribution with a step function with the same accuracy. In the Appendix, we present in detail our algorithm to approximate the ruin probability with guaranteed error bound $\delta$ by approximating the claim size distribution with accuracy of at most $\epsilon =\delta (1-\rho)/\rho$, a result which is a consequence of Proposition~\ref{ruinbound}. The exact relation between the number of phases, the accuracy $\epsilon$ and the bound $\delta$ is given also in the Appendix.

\section{Heavy-traffic and heavy-tail approximations}\label{heavy-traffic and heavy-tail approximations}

In this section, we present the heavy traffic \cite{kingman62} and the heavy tail approximations \cite{bahr75,borovkov92,embrechts82,pakes75}, which are most often used for the evaluation of the ruin probability. We first start with the heavy traffic approximation.

If the claim size distribution $B$ has a finite second moment, then as $\rho \rightarrow 1$, $M$, which was defined in Section~\ref{error bound}, converges to an exponential random variable with mean $\e M$; i.e. $Exp(1/\e M)$. This result is known as the heavy traffic approximation \cite{kalashnikov-GSBREA}. In other words,
\begin{equation}\label{heavytrafficapproximation1}
  \df{\psi} \approx \df{\psi_h} := e^{-u/\e M},
\end{equation}
where $\e M= \frac{\rho \e U^2}{2(1-\rho)\e U}$. Although the heavy traffic approximation is given through a simple exponential, its biggest drawback is that it requires the first two moments of the claim size distribution to be finite, which is not always the case for heavy-tailed distributions, e.g. the Pareto.

Equation \eqref{polllaczek khinchine formula} shows that $M$ can be written as a geometric random sum with terms distributed according to $B_0$. Bounds for exponential approximations of geometric convolutions have been obtained by Brown \cite{brown90}. Thus, we can acquire a bound for the ruin probability by applying Theorem 2.1. of \cite{brown90}, which states that the sup norm distance\footnote{The supnorm distance between two variables is actually the supnorm distance between their distributions.} between $M$ and an exponential random variable with the same mean, namely $Exp(1/\e M)$, is
\begin{equation}\label{heavy traffic bound Brown}
  \distance{M,Exp(1/\e M)}=(1-\rho)\max(2\gamma,\gamma/\rho) =  \begin{cases}
                                                  2 (1-\rho) \gamma,    &\text{if $\rho \geq \frac{1}{2}$}\\
                                                  (1-\rho)\gamma /\rho, &\text{if $0<\rho<\frac{1}{2}$},
                                            \end{cases}
\end{equation}
where $\gamma= \frac{2\e{U^3}\e U}{3 (\e{U^2})^2}$. Thus, a finite {\it third} moment is required for the claim sizes in order to guarantee a bound for the heavy traffic approximation.

When the claim sizes belong to the subexponential class of distributions \cite{teugels75}, e.g. Weibull, lognormal, Pareto, etc., the heavy tail approximation can also be used. For $u\rightarrow \infty$, the heavy tail approximation is defined as
\begin{equation*}\label{heavytailapproximation}
  \df{\psi} \approx \df{\psi_t} := \frac{\rho}{1-\rho} \disf{\com{B}_0}.
\end{equation*}
This approximation is also given by a simple formula, which this time requires only the first moment of the claim size distribution to be finite. Its drawback though is that for values of $\rho$ close to 1, or equivalently in the heavy traffic regime, the heavy tail approximation is useful only for extremely big values of $u$. For the heavy traffic setting, there exists a comparative analysis between the heavy traffic and the heavy tail approximations \cite{olvera11} in which the point at which the heavy tail approximation becomes more suitable than the heavy traffic is examined. 

In the following section, we compare the accuracy of the spectral approximation to the accuracy of the heavy traffic and the heavy tail approximations. An interesting observation with respect to the spectral approximation is that, since it decays exponentially, it converges faster to zero than any heavy-tailed distribution. Thus, at the tail the spectral approximation is expected to underestimate the ruin probability. But an overestimation of the ruin probability for small values of $u$, compensates for the underestimation at the tail, as it will be apparent in Section~\ref{experiments}.

\section{Numerical results}\label{numerical results}

In this section we implement our algorithm in order to check the accuracy of the spectral approximation. We test the spectral approximation in 3 different classes of c.m. heavy-tailed distributions: a class of long-tail distributions introduced in \cite{abate99a}, the Weibull distribution and the Pareto distribution.

\subsection{Test distributions}\label{test distributions}
First we present the three test distributions, and thereafter we do a series of experiments to compare the accuracy of the spectral approximation with the accuracy of heavy tail approximation and when applicable with the heavy traffic approximation too.

\subsubsection{Abate-Whitt distribution}\label{Abate-Whitt}
Consider a claim size probability density function \df{b}{} with Laplace transform
\begin{equation*}
  \disf[s]{\lt{b}} =1-\frac{s}{(\mu + \sqrt{s})(1 + \sqrt{s})},
\end{equation*}
which has mean $\mu^{-1}$ and all higher moments infinite. The parameter $\mu$ of the pdf $b$ can range over the positive values. This pdf was introduced in \cite{abate99a}, where it was also proven that the explicit formula for the ruin probability of the compound Poisson model with arrival rate for claims $\lambda$ and $\rho=\lambda /\mu <1$ is
\begin{align*}
  \df{\psi} &= \pr(M>u) = \frac{\rho}{v_1-v_2}\left(v_1\zeta(v^2_2 u)-v_2\zeta(v_1^2 u)\right),
\intertext{where}
       \df{\zeta} & \equiv e^u \erfct{\sqrt u},
\intertext{and}
  v_{1,2}     &=\frac{1+\mu}{2}\pm \sqrt {\left(\frac{1+\mu}{2}\right)^2-(1-\rho)\mu}.
\end{align*}
The existence of an exact formula for the ruin probability, makes this distribution very interesting because we can compare the spectral approximation with the exact ruin probability and not with the outcome of a simulation.

For this model we have that the ccdf of the the claim size distribution is given by the formula
\begin{align*}
  \disf{\com{B}}  &= \left(\frac{1}{1-\mu}\right)\left(\df{\zeta}-\mu\zeta(\mu^2 u)\right).
\intertext{With simple calculations we can verify that $\df{\zeta}$ is c.m. since it can be written as a mixture of exponentials}
       \df{\zeta} & =e^u \erfct{\sqrt u}
                    \stackrel{z=x^2}{=} \frac{2 e^u}{\sqrt \pi} \int_u^{+\infty}\frac{e^{-z}}{2\sqrt z}dz\\
                  & =\frac{1}{\sqrt \pi}\int_u^{+\infty}\frac{e^{-(z-u)}}{\sqrt z}dz
                    \stackrel{t=z-u}{=}\frac{1}{\sqrt \pi}\int_0^{+\infty}\frac{e^{-t}}{\sqrt {t+u}}dt\\
                  & =\frac{1}{\sqrt \pi}\int_0^{+\infty}\frac{e^{-t}}{\sqrt u} \left(\frac{u}{t+u}\right)^{\frac12}dt
                    =\frac{1}{\sqrt \pi}\int_0^{+\infty}\frac{e^{-t}}{\sqrt u}
                    \left(\frac{1}{\sqrt \pi}\int_0^{+\infty} \frac{\sqrt u}{\sqrt x}e^{-(u+t)y}dy\right)dt\\
                  &=\frac{1}{\pi}\int_0^{+\infty}\frac{e^{-uy}}{\sqrt y}
                    \underbrace{(\int_0^{+\infty}e^{-(y+1)t}dt)}_{\frac{1}{y+1}}dy
                   =\int_0^{+\infty}y e^{-u y} \frac{1}{\pi y^{3/2}(y+1)}dy.
\intertext{The ccdf of the claim sizes is also c.m. That is,}
  \disf{\com{B}}  &= \left(\frac{1}{1-\mu}\right)\left(\df{\zeta}-\mu\zeta(\mu^2 u)\right)\\
                  &= \frac{1}{1-\mu} \int_0^{+\infty} y e^{-u y} \left[\frac{1}{\pi y^{3/2}(y+1)}-\frac{\mu^2}{\pi y^{3/2}(y+\mu^2)}\right]dy\\
                  &=\int_0^{+\infty} e^{-u y} \frac{\sqrt{y}(1+\mu)}{\pi (y+1)(y+\mu^2)}dy.
\end{align*}

Note that for the heavy traffic approximation a finite second moment is required, which does not hold for this case. Therefore, for this distribution the heavy traffic approximation for the ruin probability cannot be evaluated. As a result, we compare the spectral approximation only with the heavy tail approximation.

\subsubsection{Weibull}\label{Weibull distribution}
The ccdf of the Weibull$(c,a)$ distribution with $c$ and $a$ the positive shape and scale parameters respectively is given by $\disf{\com{B}}=e^{-(u/a)^c}$. It can be verified \cite{jewell82} that the ccdf of the Weibull$(0.5,a)$ distribution with fixed shape parameter $\frac{1}{2}$ arises as a mixture of exponentials, where the mixing measure (measure of the spectral function) $G$ is given by
\begin{equation*}
  d\disf[y]{G} = \frac{a e^{-a^{2/4y}}}{2\sqrt{\pi y^3}}dy.
\end{equation*}

For this case we do not have an explicit formula for the ruin probability, thus we compare the spectral approximation to simulation results. Since the second moment of Weibull$(c,a)$ is finite, namely $\e B^2 = 24 a^2$, we can compare the spectral approximation with the heavy traffic approximation as well, contrary to the Abate-Whitt distribution, where only comparisons with the heavy tail approximation were possible.

\subsubsection{Pareto}\label{Pareto distribution}
The third test function we use is the Pareto$(a,b)$ distribution with shape parameter $a>0$ and scale parameter $b>0$. The Pareto$(a,b)$ distribution with pdf $\df{b}{} =\frac{ab}{(1+b u)^{a+1}}, u>0$ is completely monotone. Its ccdf $\disf{\com{B}}=(1+b u)^{-a}$ can be written as a mixture of exponentials in the form
\begin{equation*}
  (1+b u)^{-a} =\int_0^{+\infty}e^{-y u} e^{-y/b}\frac{\left(\frac{y}{b}\right)^{a-1}}{b\Gamma(a)} dy.
\end{equation*}
Also for this distribution the ruin probability does not exist in closed form. Therefore we compare our approximation for this case to simulation results.

It is known that the $n$th moment of the Pareto distribution exists if and only if the shape parameter is greater than $n$. Since it would be interesting to compare the spectral approximation, not only with the heavy tail one, but with the heavy traffic too, it is necessary to have a finite second moment for the claim sizes. Moreover, as stated in Section~\ref{heavy-traffic and heavy-tail approximations}, a bound for the heavy traffic approximation is guaranteed as long as the third moment of the distribution is finite. For these reasons, if we want to evaluate the heavy traffic approximation with a guaranteed bound for the Pareto$(a,b)$, the shape parameter $a$ must be chosen to be greater than 3.

\subsection{Numerical results}\label{experiments}
The goal of this section is to implement our algorithm to check the accuracy of the spectral approximation and the tightness of its accompanying bound, which is given in Proposition~\ref{ruinbound}. More precisely, we answer the following questions.

Since the only restriction we have for the parameters of the three test distributions is that the shape parameter of the Pareto$(a,b)$ must be greater than 3, we randomly select the parameters and thus we deal with the Abate-Whitt distribution with $\mu =2$, the Weibull$(0.5,3)$ distribution and the Pareto$(4,3)$ distribution.

\begin{enumerate}
  \item \textsc{Impact of phases}. The bound of the spectral approximation is conversely proportional to the number of phases of the hyperexponential with which we approximate the excess claim size distribution (see Appendix). So, for a fixed claim rate $\rho$, the bound becomes tighter when the number of phases increases. Does this also mean that the spectral approximation becomes more accurate as the number of phases increases?\label{Q.increase the number of phases}

      \textsc{Experiment}: We fix $\rho$ and we compare three different spectral approximations with number of phases 10, 20 and 100 respectively, with the exact value of the ruin probability. For the Abate-Whitt distribution, we present the exact ruin probability with the three approximations in one graph; see Figure~\ref{figure:Abate increase phases}. For the Weibull and the Pareto distributions we compare the three approximations to the exact ruin probability that we obtain through simulation and display our results in Tables~\ref{table:Weibull for different number of phases} and \ref{table:Pareto for different number of phases}. As for all different values of $\rho$ we get a similar results, we present our findings only for $\rho = 0.7$.

      \textsc{Answer}: The conclusion is that, while the number of phases increases, a more accurate spectral approximation is achieved. This result is in line with our expectations, and we can safely conclude that for a fixed claim rate $\rho$ more phases lead to a more accurate spectral approximation.

  \item \textsc{Quality of the bound}. Is the bound strict or pessimistic? How far is the bound from the real error of the spectral approximation? \label{Q.quality of error bound}

      \textsc{Experiment}: We fix the bound of the spectral approximation to be equal to $\delta = 0.02$, and we evaluate the error functions (in absolute values) for the spectral approximation when the claim rate $\rho$ takes the values 0.1, 0.5 and 0.9. For these three cases we need 5, 49, and 449 phases respectively for the spectral approximation. We compare the guaranteed bound with the exact maximum error that is achieved; see Figures~\ref{figure:Abate all error functions} to \ref{figure:Pareto all error functions}.

      Also, for various combinations of number of phases and claim rate $\rho$, we calculate the ratios between the predicted bound of the spectral approximation and the achieved maximum error; see Table~\ref{table:Abate with ratios}. We set out this experiment only for the Abate-Whitt distribution, because the existence of the exact ruin probability gives more accurate results.

      \textsc{Answer}: An interesting observation that arises from Figure~\ref{figure:Abate all error functions} is that the achieved maximum error of the spectral approximation seems to be almost half of the guaranteed bound. In order to verify that the bound is twice as big as the achieved maximum error we look at Table~\ref{table:Abate with ratios}.

      We first read the table horizontally, namely we fix the claim rate $\rho$. We observe that while we let the number of phases increase, the ratio between the predicted bound and the real maximum error becomes smaller and converges to 2. As it was mentioned earlier, the spectral approximation becomes more accurate when we increase the number of phases. Therefore, we conclude that the bound becomes tighter when for a fixed $\rho$ we increase the number of phases.

      We read now the table vertically, namely we fix the number of phases and we let the claim rate $\rho$ increase. We observe that while we let $\rho$ increase, both the predicted bound and the maximum error increase. Since the ratios between the bound and the maximum error increase too, we can conclude that the bound becomes less tight when the claim rate increases.

      However, from Figures~\ref{figure:Weibull all error functions} and \ref{figure:Pareto all error functions}, we see that the achieved maximum error is not only 2 times smaller than the guaranteed bound but 4 times smaller! Gathering all the above together, we can conclude that the bound seems to be at least twice as big as the the achieved maximum error of the spectral approximation.

  \item \textsc{Comparison of Spectral, Heavy tail, Heavy traffic approximations}. The accuracy of the spectral approximation can be predetermined through its bound. For a fixed range of $u$, which of the three approximations -- spectral, heavy tail and heavy traffic (when applicable) -- is better than the others as $\rho \rightarrow1$ or $\rho\rightarrow0$, when the bound predicts accuracy of at most $\delta$ for the spectral approximation? \label{Q.comparison approximations}

      \textsc{Experiment}: We fix the bound of the spectral approximation to be equal to $\delta = 0.02$, and for $\rho=0.1,0.5$ and 0.9 we compare the spectral (with 5, 49 and 449 phases respectively), the heavy tail and the heavy traffic (when applicable) approximations. We present the distributions in a graph, where the displayed range of $u$ is such that $\df[u]{\psi}>\delta$, because after this point the error is smaller than $\delta$. The level $\delta$ is denoted on the graphs with a dashed horizontal line; see Figures ~\ref{figure:Abate-Whitt 0.1, with accuracy 0.02} to \ref{figure:Pareto 0.9, with accuracy 0.02}.

      \textsc{Answer}: We observe that the spectral approximation behaves nicely for all values of $u$. For small values of $u$, the spectral approximation is more accurate than the heavy tail approximation, where the second fails to provide us with a good estimation of the ruin probability, especially when $\rho \rightarrow 1$.

      On the other hand, the heavy tail approximation is slightly more accurate than the spectral approximation at the tail. Although we cannot give an estimation for the point $u^*$ at which the heavy tail approximation becomes more suitable than the spectral approximation, we observe that this point takes greater values as $\rho$ increases and it sometimes can be extremely big; i.e. see Figure~\ref{figure:Abate-Whitt 0.9, with accuracy 0.02}.

      Furthermore, according to our expectations, the spectral approximation overestimates the ruin probability for small values of $u$ (this is more clear for small values of $\rho$) and underestimates it for large values of $u$. In all cases, the heavy traffic approximation is worse than the other two, since it exhibits a sharper behavior than the spectral approximation. Namely, for small values of $u$ it overestimates the ruin probability more than the spectral approximation, and for large values of $u$ it underestimates the ruin probability more than the spectral approximation. Note also that, at the tail, the spectral approximation and the heavy traffic approximation are almost identical, which can be explained by the fact that both of them have an exponential decay.

  \item \textsc{Comparison between Spectral and Heavy traffic bounds}. For the Weibull and the Pareto distributions, the heavy traffic approximation can be evaluated and it also has a guaranteed bound \cite{brown90}. So, is there a rule of thumb to help us choose between the spectral and the heavy traffic approximation, when they both guarantee the same bound?

      \textsc{Experiment}: For various values of $\rho$, we compare the spectral approximation with the heavy traffic approximation when they both guarantee the same bound. More precisely, we fix $\rho$ and determine the number of phases $k^*$ of the spectral approximation for which both approximations guarantee the same bound. We calculate the two approximations and evaluate their maximum errors. We present our findings in a table, only for some values of $\rho$ that the heavy traffic bound has a meaning, namely when it is smaller than 1; see Tables~\ref{table:Weibull comparison spectral with heavy traffic} and \ref{table:Pareto comparison spectral with heavy traffic}.

      We can easily verify that for the Pareto$(a,b)$ distribution, the heavy traffic bound depends on the shape parameter $a$, since $\gamma = \frac{a-2}{a-3}$. An interesting experiment that arises from this observation is to check whether we have a clearer picture on which of the spectral and heavy traffic approximations is the best in terms of accuracy, if we choose $a$ big enough such that $\gamma \rightarrow 1$, namely if we make the heavy traffic bound tighter (for Pareto$(4,3)$, $\gamma =2$). For this reason, we repeat our last experiment for Pareto$(15.6,2.7)$, which has $\gamma =1.079$.

      \textsc{Answer}: From Table~\ref{table:Weibull comparison spectral with heavy traffic}, which gives the results for Weibull$(0.5,3)$, we see that whenever the bounds are equal, the spectral approximation is more accurate than the heavy traffic approximation for all number of phases  greater or equal than $k^*$. On the other hand, from Table~\ref{table:Weibull comparison spectral with heavy traffic}, which gives the results for Pareto$(4,3)$, we get a different picture. The conclusion that we draw from this table is that for a small number of phases (relatively smaller than 20), the heavy traffic approximation is better while for a number of phases greater than 20 the conclusion reverses.

      For Pareto$(15.6,2.7)$, more phases were needed in the corresponding spectral approximation for the same values of $\rho$, because the heavy traffic bound is now tighter. The picture from Table~\ref{table:Pareto comparison spectral with heavy traffic 2} is not that clear. More precisely, even when the number of phases becomes relatively big we cannot draw a safe conclusion that the spectral approximation is better than the heavy traffic approximation.
\end{enumerate}

 At this point it is interesting to observe the following. The heavy traffic approximation as presented in Section~\ref{heavy-traffic and heavy-tail approximations} has no atoms. It is known \cite{asmussen-RP} that the ruin probability has an atom of mass $\rho$ at 0. Thus, the heavy traffic approximation is not very accurate for small values of $u$, especially when $\rho$ takes relatively small values. For this reason, a more suitable heavy traffic approximation ($\psi_h$) for our comparisons for all values of $\rho$ seems to be
        \begin{equation}\label{heavytrafficapproximation}
          \df{\psi} \approx \df{\psi_h} := \rho e^{-\rho u/\e M},
        \end{equation}
for which is easy to verify that it also has mean equal to $\e M$ and an atom of mass $\rho$ at 0. Since we used a different heavy traffic approximation in all of our experiments than the one Brown \cite{brown90} compares the ruin probability with, we extended Brown's bound, given in \eqref{heavy traffic bound Brown}, to this situation. Applying the triangular inequality to the sup norm distance we get
\begin{equation*}
    \distance{\psi,\psi_h} \leq \distance{\psi,Exp(1/\e M)} + \distance{Exp(1/\e M),\psi_h}.
\end{equation*}
It is easy to verify that $\distance{Exp(1/\e M),\psi_h} = 1-\rho$, so the sup norm distance between the ruin probability and the heavy traffic approximation we use for comparisons is
\begin{equation}\label{heavy traffic bound}
    \distance{\psi , \psi_h} \leq (1-\rho)\max(2\gamma,\gamma/\rho) + 1-\rho = (1-\rho) \cdot \begin{cases}
                                                                             2  \gamma + 1,    &\text{if $\rho \geq \frac{1}{2}$}\\
                                                                             \gamma/\rho + 1, &\text{if $0<\rho<\frac{1}{2}$},
                                                                          \end{cases}
\end{equation}
where $\gamma= \frac{2\e{U^3}\e U}{3 (\e{U^2})^2}$. When we referred to the heavy traffic approximation and its accompanying bound, in all of our experiments we meant those given from \eqref{heavytrafficapproximation} and \eqref{heavy traffic bound}, respectively.

\section{Conclusions}\label{conclusions}
In this paper we addressed the problem of how many phases are needed to approximate a heavy-tailed distribution with a phase-type distribution in such a way that one can obtain a guaranteed bound on the approximation of the ruin probability (see Appendix). In doing so, we developed an explicit bound using the geometric random sum representation, which was combined with a spectral approximation of the excess claim size distribution.

The conclusions that we can draw, both for the spectral approximation and its bound, can be summarized as follows:
\begin{itemize}
  \item The spectral approximation provides a good fit for all values of $u$, especially for the small ones, where the heavy traffic and heavy tail approximations fail. Also, for small values of $u$ the spectral approximation exhibits a behavior of overestimating the ruin probability, while for larger values of $u$ we have an underestimation of the ruin probability by the spectral approximation. Finally, for a fixed claim rate $\rho$, the more the phases we have for the approximate hyperexponential of the excess claim size distribution, the more accurate spectral approximation we achieve.
  \item The spectral bound, guaranteed by Proposition~\ref{ruinbound}, becomes tighter when for a fixed claim rate $\rho$ the number of phases is increased, while it becomes less tight when for a fixed number of phases the claim rate increases. Moreover, the bound seems to be at least twice as big as the achieved maximum error of the spectral approximation. But, based only on the numerical examples we performed, we cannot conclude that this is the general rule.
  \item Based on existing analytical results and extensive experiments it is hard to draw a definitive conclusion on which approximation should be preferred: the heavy traffic approximation or the spectral approximation. We believe that obtaining more mathematical as well as experimental insights in this  problem is an important topic for future research.
\end{itemize}

To sum up, the spectral approximation provides a good fit for all values of $u$ and has a guaranteed accuracy, while it requires only a finite mean for the claim sizes.

\section*{Acknowledgments}\label{acknowledgments}
The work of Maria Vlasiou and Eleni Vatamidou is supported by Netherlands Organisation for Scientific Research (NWO) through project number 613.001.006. The work of Bert Zwart is supported by an NWO VIDI grant and an IBM faculty award.

\bibliographystyle{amsplain}
\bibliography{Vatamidou_Eleni}

\begin{appendix}
  \appendix

\section{Algorithm}

We consider the compound Poisson model \cite{asmussen-RP} with a completely monotone claim size distribution $B$ with a finite mean. We denote by $B_0$ the excess claim size distribution and by $G_0$ the latter's spectral function (strictly increasing distribution). We develop an algorithm to evaluate the ruin probability by approximating the excess claim size distribution with a hyperexponential one, where all phases have equal weights.

Before we present the spectral approximation algorithm, it is necessary to give an important property on which the Laplace inversion of the ruin probability will be based on. The Pollaczeck-Khinchine formula \eqref{laplace transform of polllaczek khinchine formula} can be written equivalently in the form
\begin{equation}
  \e{[e^{-sM}]} = 1-\rho + \rho\frac{(1-\rho)\disf[s]{\lt{B}_0}}{1-\rho \disf[s]{\lt{B}_0}} = 1-\rho+\rho \disf[s]{\lt{{M}}_+},
\end{equation}
where $M_+ \stackrel{d}{=} M|M>0$. We have the following lemma.

\begin{lemma}\label{real roots Cohen}
 If $B_0$ follows a hyperexponential distribution with $k$-phases then $M_+$ follows a hyperexponential distribution with $k$-phases as well (with different exponential rates and weights from the first one). In other words, $\disf[s]{\lt{{M}}_+}$ can be written in the form $\sum_{i=1}^k R_i\frac{\eta_i}{\eta_i+s}$ for some $R_i,\eta_i, i=1,...,k$.
\end{lemma}
\begin{proof}
In \cite{cohen-SSQ} it was proven that $M_+$ follows a hyperexponential distribution in the $G/K_k/1$ queue. Since the $M/H_k/1$ is a special case of the $G/K_k/1$ queue, the result holds here too.
\end{proof}

After this, we present our algorithm

{\bf begin algorithm}
  \begin{enumerate}
    \item Write $\disf{\com{B}}$ as a mixture of exponentials.
    \item Using Lemma~\ref{the complementaty integrated tail}, find the spectral function $\disf[y]{G_0}$ of $\disf{\com{B}_0}$.
    \item Approximate $\disf{\com{B}_0}$ by a hyperexponential distribution with $k$-phases.
        \begin{enumerate}
          \item Choose the number of phases $k$.
          \item Set the accuracy of the approximation $\epsilon = \frac{1}{k+1}$, such that $\left|\disf{B_0}-\disf{\apr{B}_0} \right| \leq \epsilon$.\label{step}
          \item Define $k$ quantiles such that $\disf[\lambda_i]{G_0} = i\epsilon$, $i=1,...,k$.\label{quantiles}
          \item  Approximate the spectral function by the step function
                \begin{equation*}
                    \disf[y]{\apr{G}_0} =\begin{cases}
                              0,             & y\in[0,\lambda_1),\\
                              \frac{i}{k},   & y\in [\lambda_i,\lambda_{i+1}),\\
                              1,             & y\geq \lambda_k.
                          \end{cases}
                \end{equation*}
          \item Find the approximation of the excess claim size distribution as $\disf[y]{\com{\apr{B}}_0}=\frac{1}{k}\sum_{i=1}^k e^{-y \lambda_i}$.
        \end{enumerate}
    \item Calculate its Laplace transform $\disf[s]{\lt{\apr{B}}_0} = \frac{1}{k}\sum_{i=1}^k\frac{\lambda_i}{\lambda_i+s}$.
    \item Choose $\rho$.
    \item Calculate $\disf[s]{\lt{\apr{{M}}}_+}$ through the formula  $\disf[s]{\lt{\apr{M}}_+} = \frac{(1-\rho)\disf[s]{\lt{\apr{B}}_0}}{1-\rho \disf[s]{\lt{\apr{B}}_0}}$.
    \item Using Lemma~\ref{real roots Cohen}, split $\disf[s]{\lt{\apr{{M}}}_+}$ into simple fractions and estimate their roots $\eta_i$. Calculate also the coefficients $R_i$.
    \item Invert the Laplace transform of $\e{[e^{-sM}]}  = 1-\rho+\rho\disf[s]{\lt{\apr{{M}}}_+}$ and find that
        $\df{\apr{\psi}} = 1-\rho +\rho\sum_{i=1}^k R_i(1-~e^{-\eta_i u})$, $u\geq0$.
    \item The accuracy for $\df{\apr{\psi}}$ is then $\delta = \epsilon \frac{\rho}{1-\rho}$.
  \end{enumerate}
{\bf end algorithm}

\begin{remark}
  At step~\eqref{step} of the algorithm we approximate the spectral function $G_0$ with a step function where the jumps occur at the quantiles $\lambda_i$ and they are all of size
  $\epsilon+\frac{\epsilon^2}{1-\epsilon}$. It can be very easily verified that, by this choice of jumps, we avoid any atoms at 0 and we still achieve $\distance{G_0,\apr{G}_0} \leq \epsilon$.
\end{remark}

\begin{remark}
  Note that the algorithm was presented under the setting that we first fix the accuracy $\epsilon$ for the approximation of the excess claim size distribution and then we evaluate the bound $\delta$ of the spectral approximation. With slight modifications, the algorithm can be presented by first fixing the desired accuracy $\delta$ for the approximation of the ruin probability. In this setting, we would have to find the number of required phases as $k = \ulcorner \frac{\rho}{(1-\rho) \delta}\urcorner -1$, where $\ulcorner x\urcorner$ is the integer which is greater than or equal to $x$ but smaller than $x+1$.
\end{remark}

\begin{remark}
  From the structure of the algorithm it is evident that we only need to write the c.m. claim size distribution as a mixture of exponentials. In comparison with the distributions used in the examples, there exist mixed distributions, such as the hyperexponential, that have a spectral function which is not strictly increasing and/or has jumps. In these cases, the algorithm cannot be applied as is and more attention needs to be paid. The problem appears at step \eqref{quantiles}, when we invert the spectral function to find the quantiles. More precisely, in a non strictly increasing spectral function we might have $\disf[x]{G_0} = i\epsilon$, for $x\in(a,b)$, with $a\neq b$, for some $i=1,...,k$. Therefore, since inversion will not give a unique value for the quantile $\lambda_i$, there must be a concrete way to define it. Also, when there are jumps, we might encounter the problem that $\disf[x]{G_0} \neq i\epsilon$ for all $x\in (0,\infty)$. In this case, $\lambda_i$ could take the value at which the jump occurred. All the above mentioned problems related to the determination of the quantiles can be overcome with small modifications to the algorithm.
\end{remark}

\section{Graphs-Tables}

\begin{figure}[ht!]
    \begin{center}
      \includegraphics[scale=0.9]{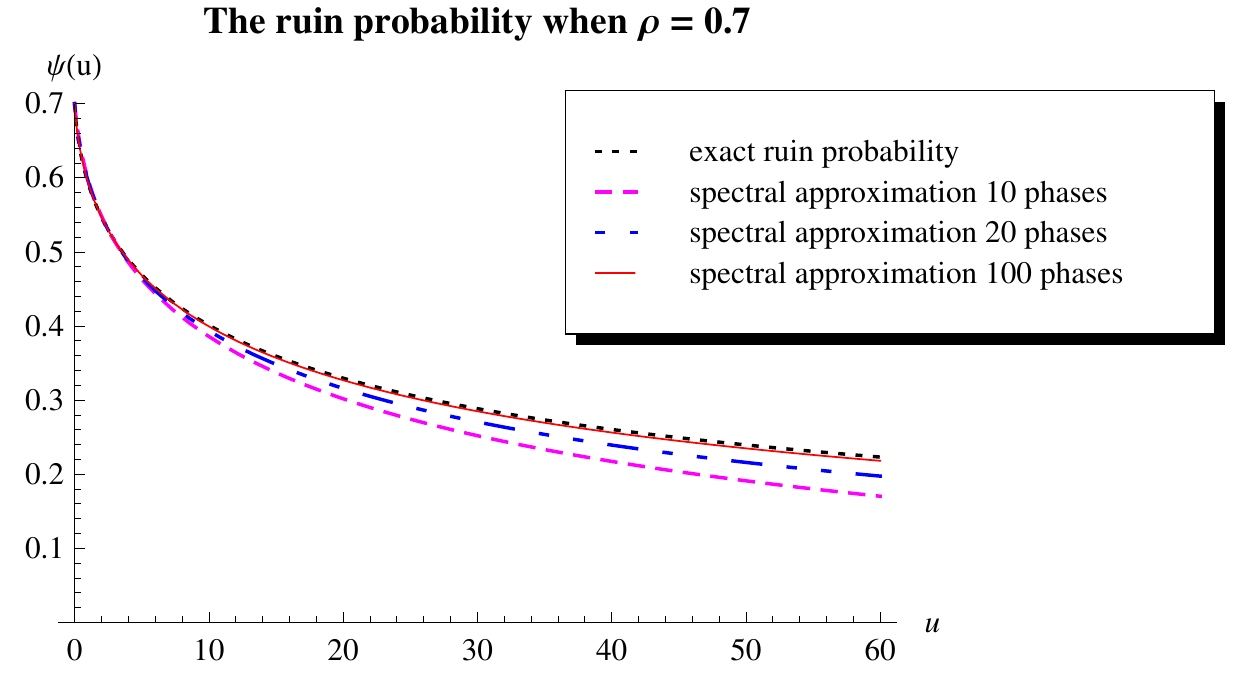}
    \end{center}
    \caption{The spectral approximation for different number of phases, when the claims follow the Abate-Whitt distribution with $\mu=2$.}\label{figure:Abate increase phases}
\end{figure}

\begin{table}[ht!]
  \begin{center}
    \begin{tabular}{|c|cccc|}\hline
      $u$    & \text{exact ruin probability}   & \text{sa 10 phases}                  & \text{sa 20 phases}                     & \text{sa 100 phases} \\ \hline
      0      & 0.70000                         & 0.70000 (\textcolor{red}{0.00000})   & 0.70000 (\textcolor{red}{0.00000})      & 0.70000 (\textcolor{red}{0.00000}) \\
      5      & 0.60745                         & 0.61023 (\textcolor{red}{0.00279})   & 0.60823 (\textcolor{red}{0.00079})      & 0.60754 (\textcolor{red}{0.00009}) \\
      10     & 0.54574                         & 0.54696 (\textcolor{red}{0.00122})   & 0.54569 (\textcolor{red}{0.00005})      & 0.54527 (\textcolor{red}{0.00047}) \\
      15     & 0.49580                         & 0.49558 (\textcolor{red}{0.00022})   & 0.49502 (\textcolor{red}{0.00078})      & 0.49485 (\textcolor{red}{0.00095}) \\
      20     & 0.45312                         & 0.45172 (\textcolor{red}{0.00139})   & 0.45181 (\textcolor{red}{0.00130})      & 0.45189 (\textcolor{red}{0.00122}) \\
      25     & 0.41603                         & 0.41334 (\textcolor{red}{0.00269})   & 0.41405 (\textcolor{red}{0.00198})      & 0.41436 (\textcolor{red}{0.00167}) \\ \hline
    \end{tabular}
  \end{center}
  \caption{The spectral approximation for different number of phases, when the claims follow the Weibull$(0.5,3)$ distributions. The numbers in the brackets correspond to the absolute error of the exact ruin probability from its respective approximations.}\label{table:Weibull for different number of phases}
\end{table}

\begin{table}[ht!]
  \begin{center}
    \begin{tabular}{|c|cccc|}\hline
      $u$     & \text{exact ruin probability}   & \text{sa 10 phases}                  & \text{sa 20 phases}                     & \text{sa 100 phases} \\ \hline
      0.00    & 0.70000                         & 0.70000 (\textcolor{red}{0.00000})   & 0.70000 (\textcolor{red}{0.00000})      & 0.70000 (\textcolor{red}{0.00000}) \\
      0.10    & 0.54805                         & 0.55012 (\textcolor{red}{0.00207})   & 0.55008 (\textcolor{red}{0.00203})      & 0.55005 (\textcolor{red}{0.00200}) \\
      0.55    & 0.23572                         & 0.22698 (\textcolor{red}{0.00873})   & 0.23218 (\textcolor{red}{0.00353})      & 0.23435 (\textcolor{red}{0.00137}) \\
      1.00    & 0.11499                         & 0.10194 (\textcolor{red}{0.01305})   & 0.10851 (\textcolor{red}{0.00648})      & 0.11146 (\textcolor{red}{0.00352}) \\
      1.45    & 0.05983                         & 0.04695 (\textcolor{red}{0.01287})   & 0.05265 (\textcolor{red}{0.00718})      & 0.05545 (\textcolor{red}{0.00437}) \\
      1.90    & 0.03215                         & 0.02187 (\textcolor{red}{0.01028})   & 0.02609 (\textcolor{red}{0.00606})      & 0.02838 (\textcolor{red}{0.00377}) \\ \hline
    \end{tabular}
  \end{center}
  \caption{The spectral approximation for different number of phases, when the claims follow the Pareto$(4,3)$ distributions. The numbers in the brackets correspond to the absolute error of the exact ruin probability from its respective approximations.}\label{table:Pareto for different number of phases}
\end{table}

\begin{figure}[ht]
\centering
    \subfigure[The claim size distribution is the Abate-Whitt distribution with $\mu=2$.]
    {
        \includegraphics[scale=0.9]{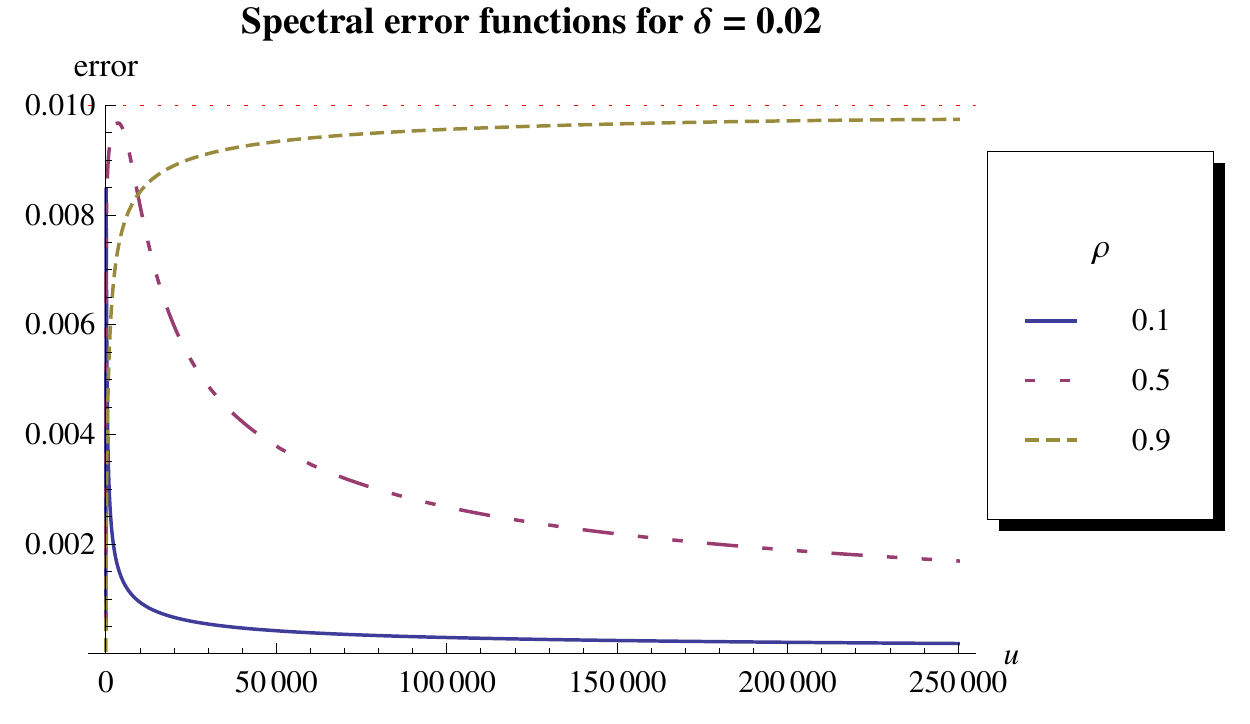}\label{figure:Abate all error functions}
    }
    \subfigure[The claim size distribution is Weibull$(0.5,3)$.]
    {
        \includegraphics[scale=0.9]{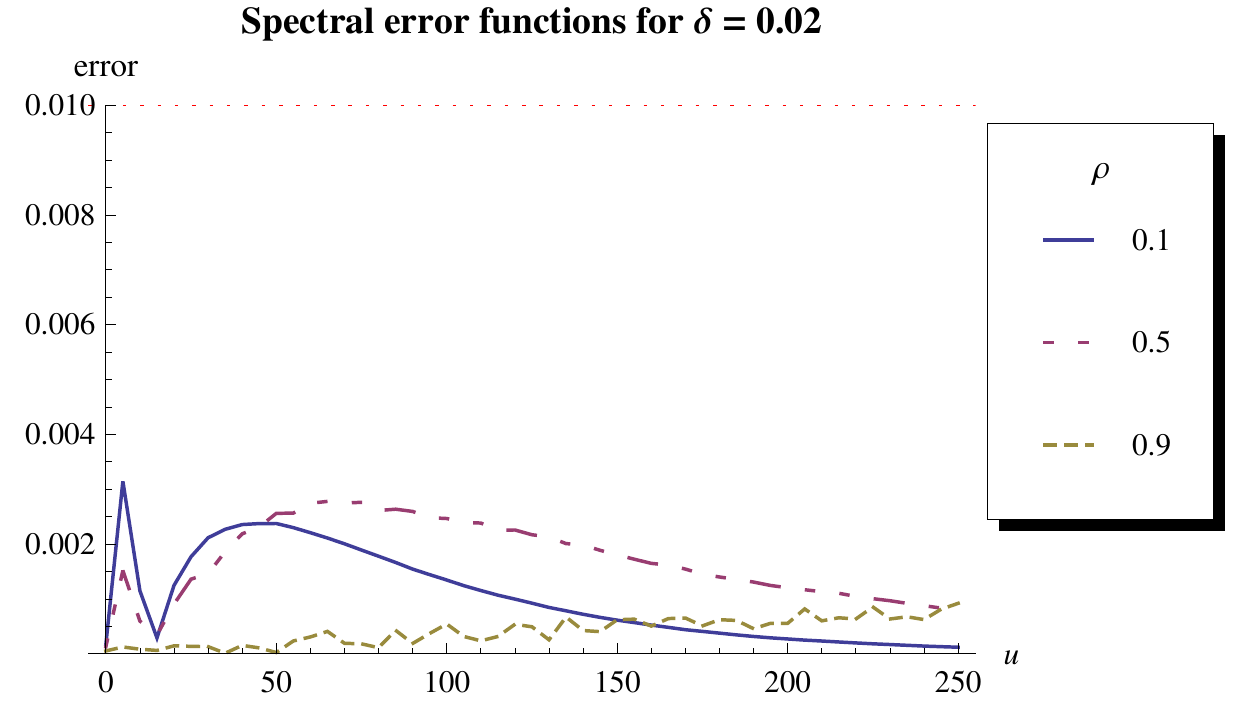}\label{figure:Weibull all error functions}
    }
    \subfigure[The claim size distribution is Pareto$(4,3)$.]
    {
        \includegraphics[scale=0.9]{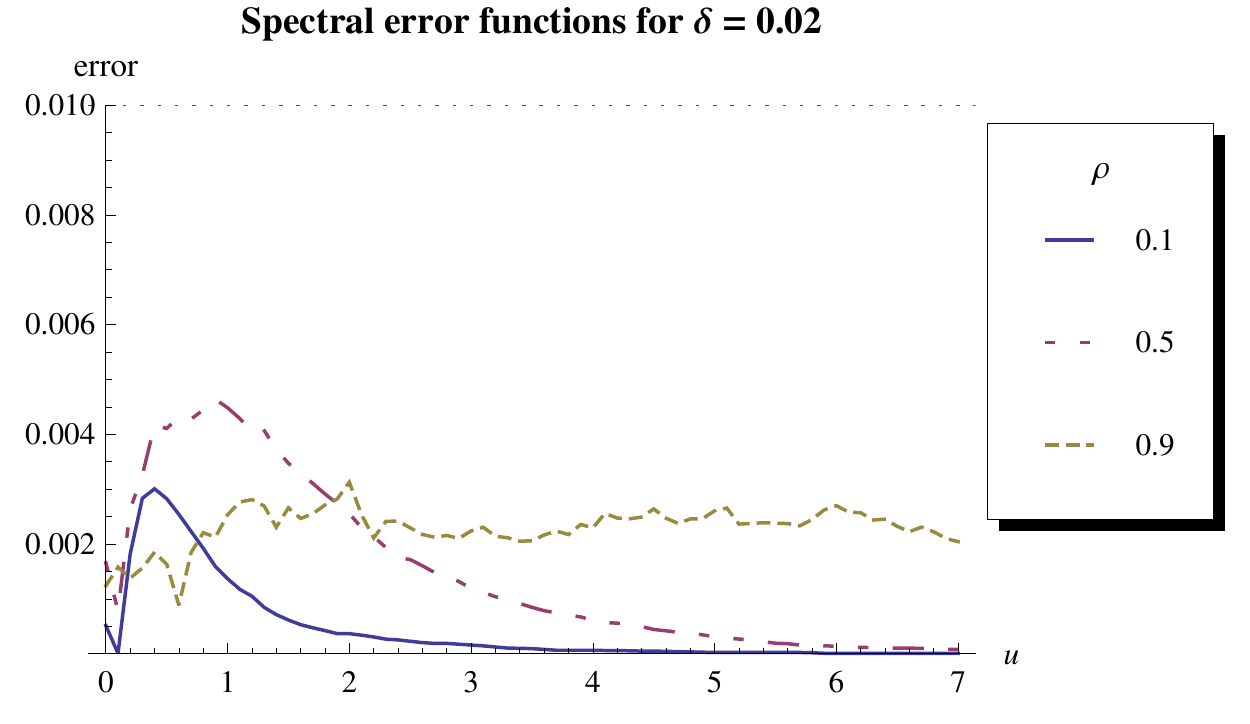}\label{figure:Pareto all error functions}
    }
\label{fig:subfigureExample}
\caption[]{Error functions for the spectral approximation with guaranteed bound $\delta =0.02$, when the claims follow each of the above distributions.}
\end{figure}

\begin{figure}[ht]
\centering
    \subfigure[]
    {
        \includegraphics[scale=0.9]{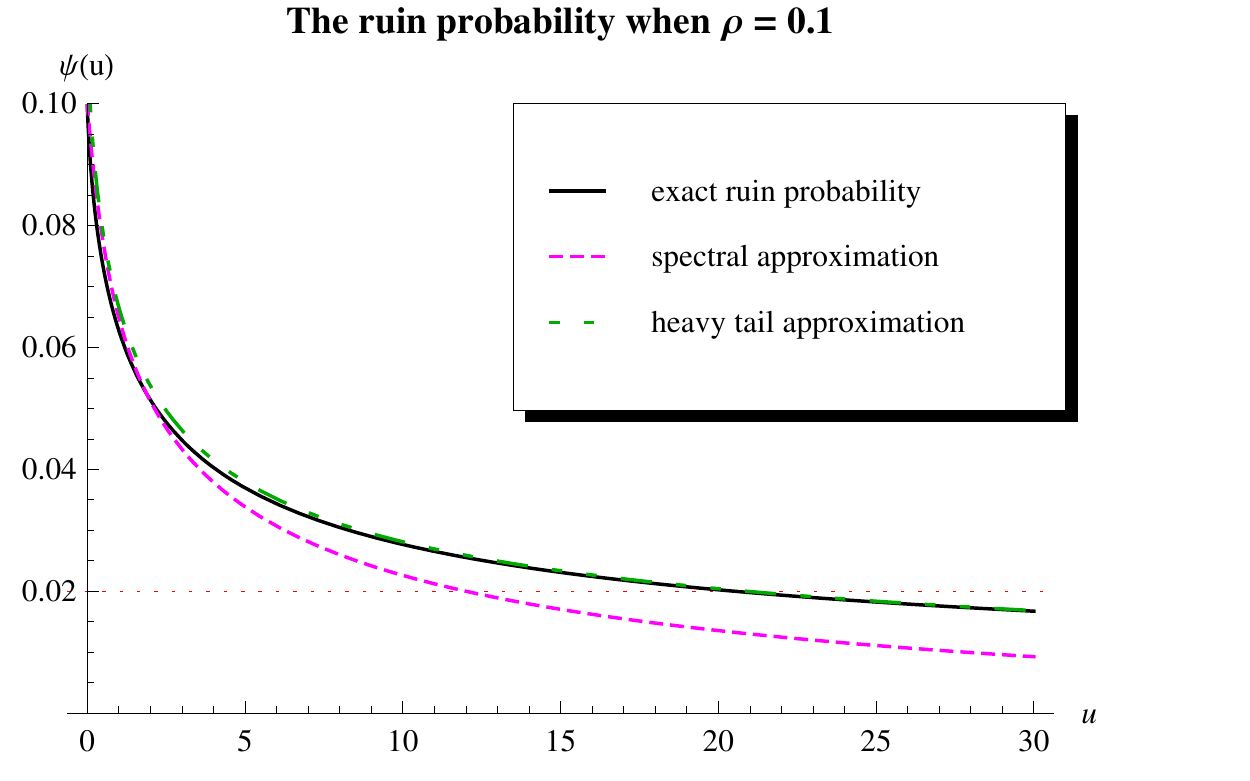}\label{figure:Abate-Whitt 0.1, with accuracy 0.02}
    }
    \subfigure[]
    {
        \includegraphics[scale=0.9]{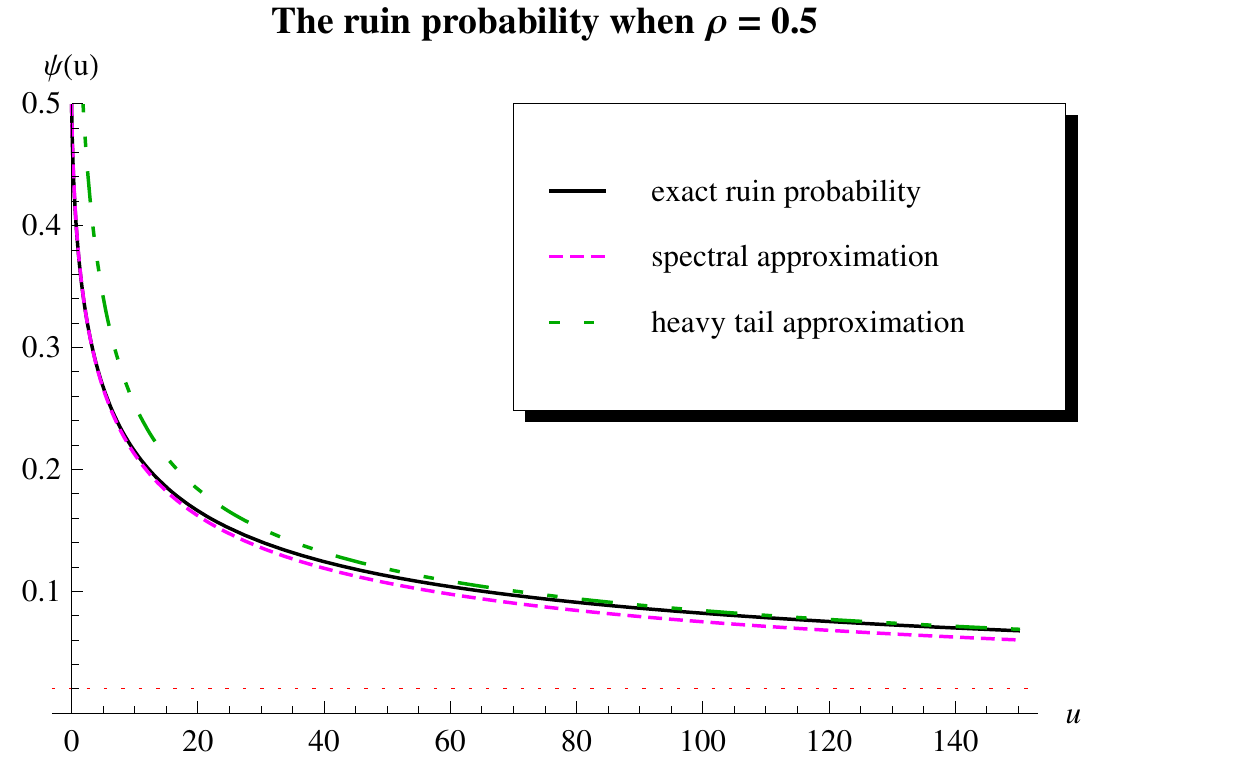}\label{figure:Abate-Whitt 0.5, with accuracy 0.02}
    }
    \subfigure[]
    {
        \includegraphics[scale=0.9]{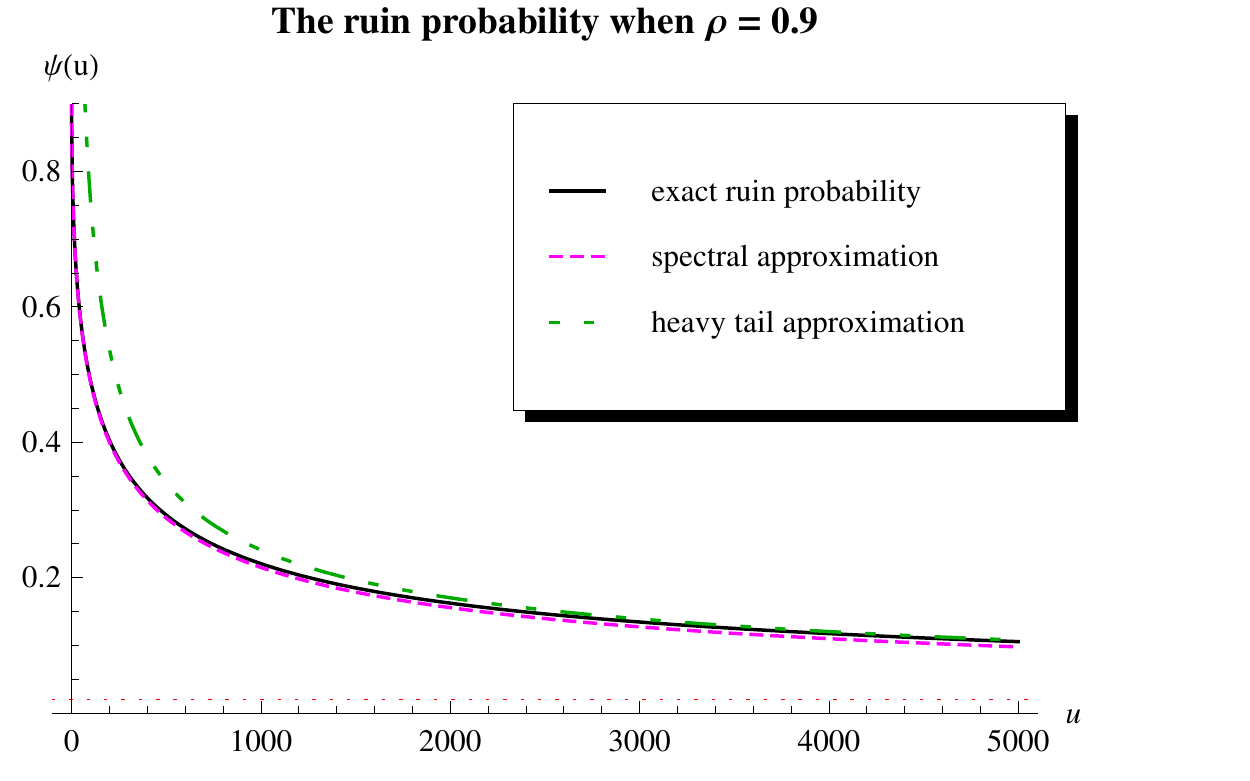}\label{figure:Abate-Whitt 0.9, with accuracy 0.02}
    }
\label{fig:subfigureExample}
\caption[]{The spectral approximation for guaranteed bound $\delta =0.02$, when the claims follow the Abate-Whitt distribution with $\mu=2$.}
\end{figure}

\begin{figure}[ht]
\centering
    \subfigure[]
    {
        \includegraphics[scale=0.9]{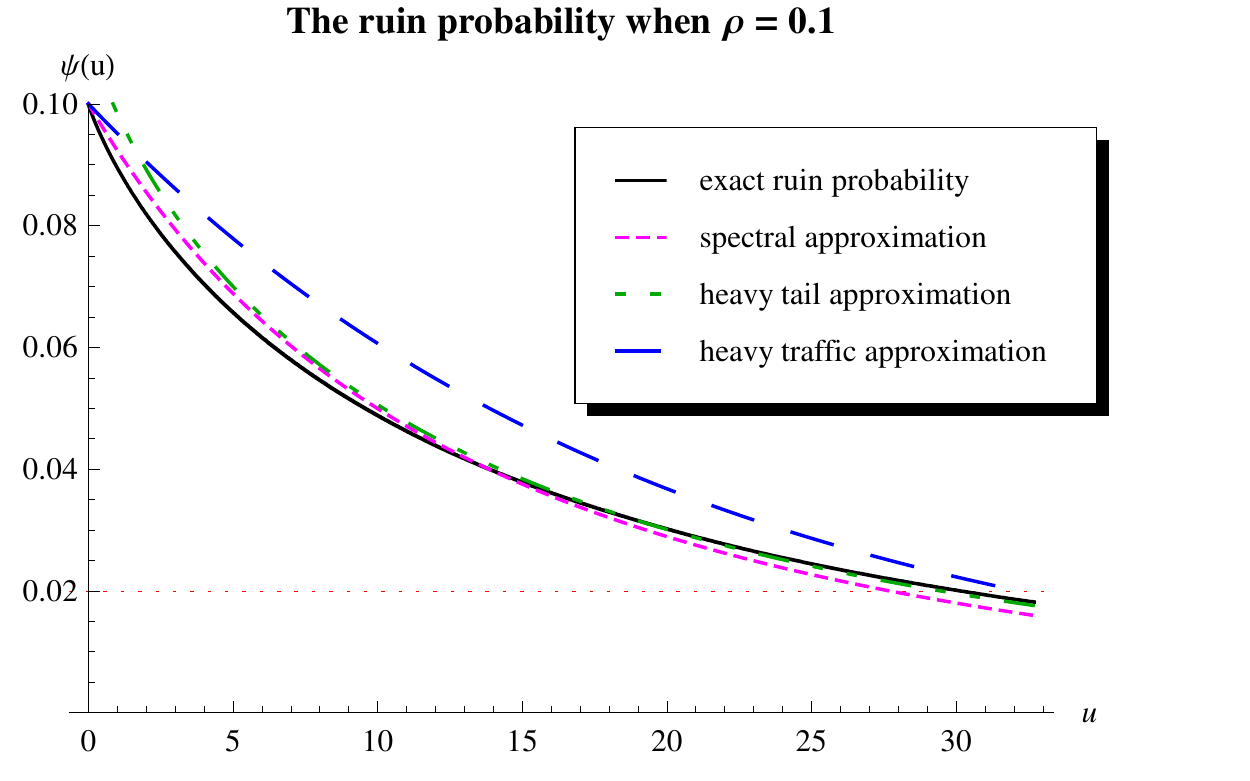}\label{figure:Weibull 0.1, with accuracy 0.02}
    }
    \subfigure[]
    {
        \includegraphics[scale=0.9]{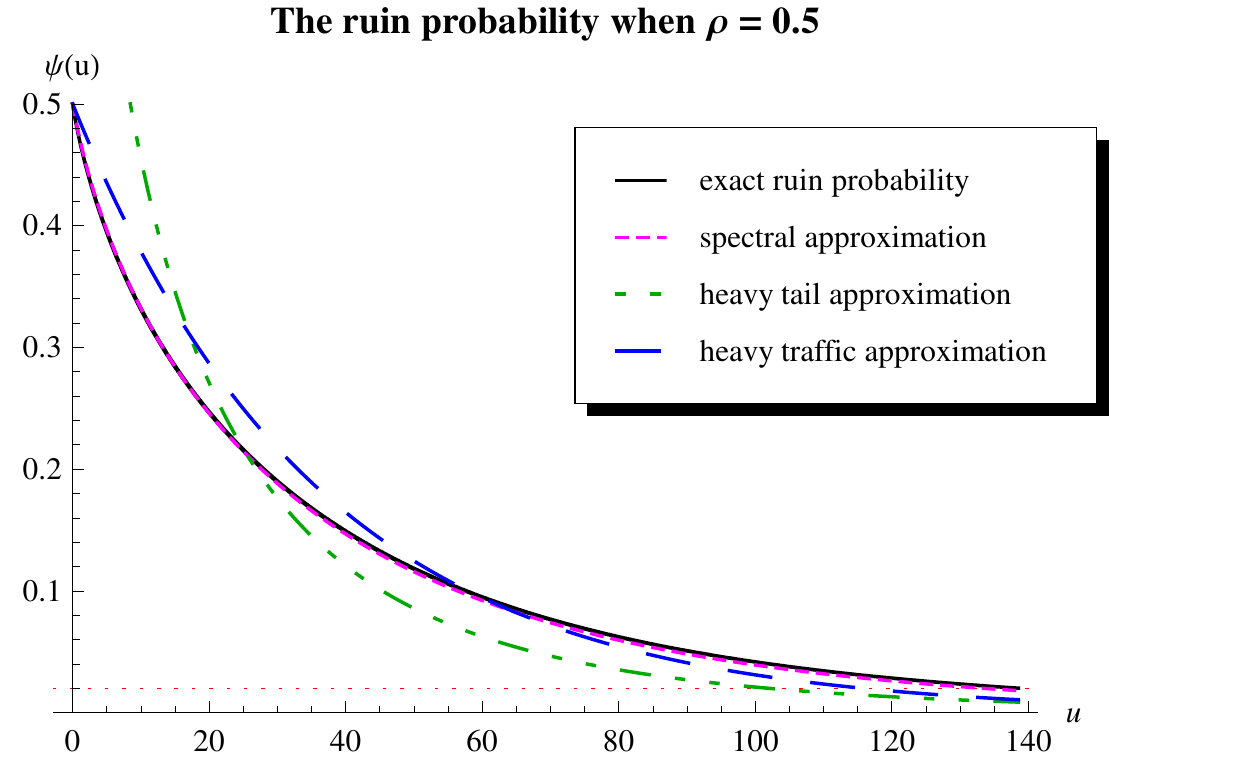}\label{figure:Weibull 0.5, with accuracy 0.02}
    }
    \subfigure[]
    {
        \includegraphics[scale=0.9]{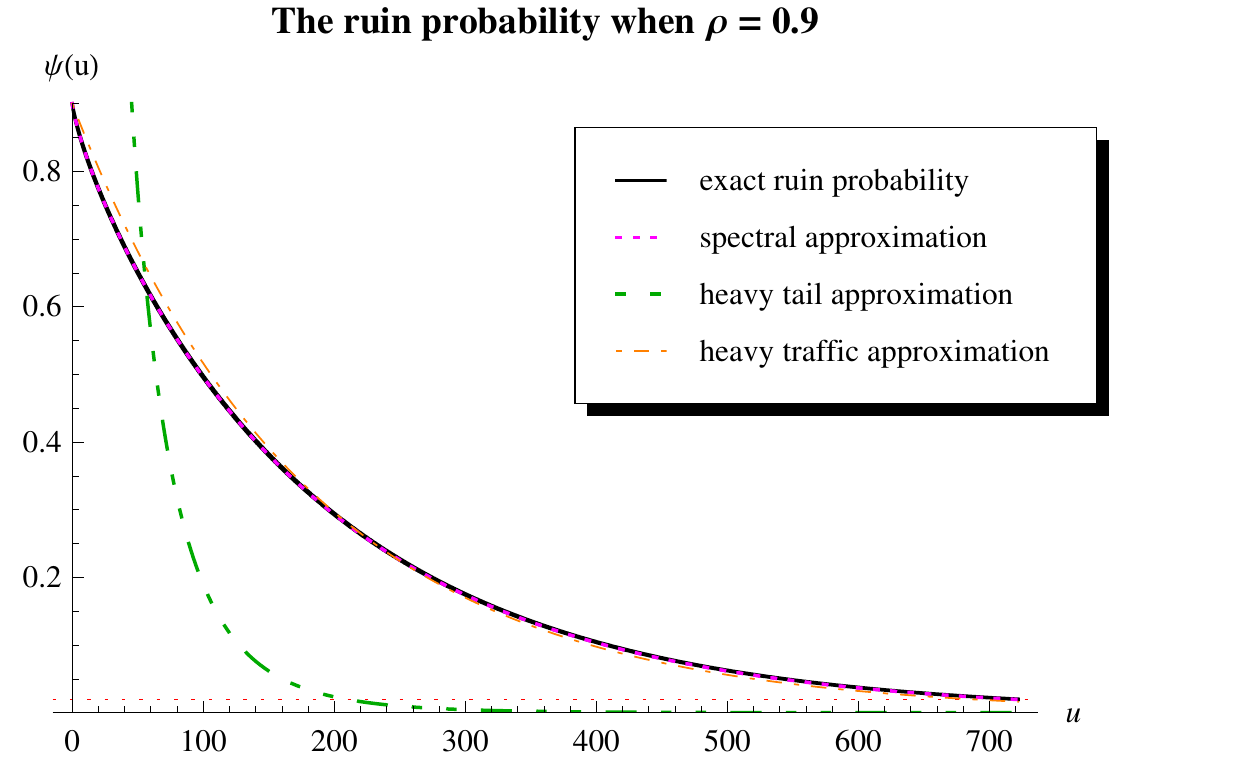}\label{figure:Weibull 0.9, with accuracy 0.02}
    }
\label{fig:subfigureExample}
\caption[]{The spectral approximation for guaranteed bound $\delta =0.02$, when the claims follow the Weibull$(0.5,3)$ distribution.}
\end{figure}

\begin{figure}[ht]
\centering
    \subfigure[]
    {
        \includegraphics[scale=0.9]{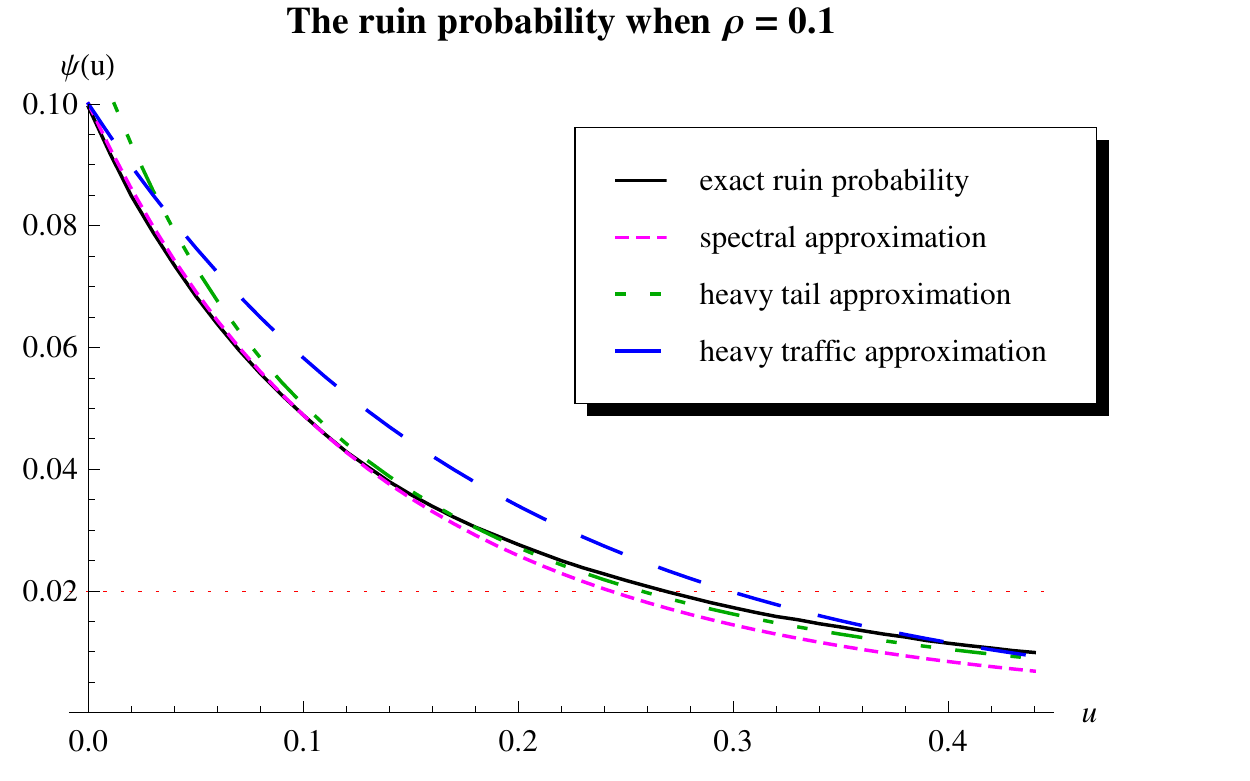}\label{figure:Pareto 0.1, with accuracy 0.02}
    }
    \subfigure[]
    {
        \includegraphics[scale=0.9]{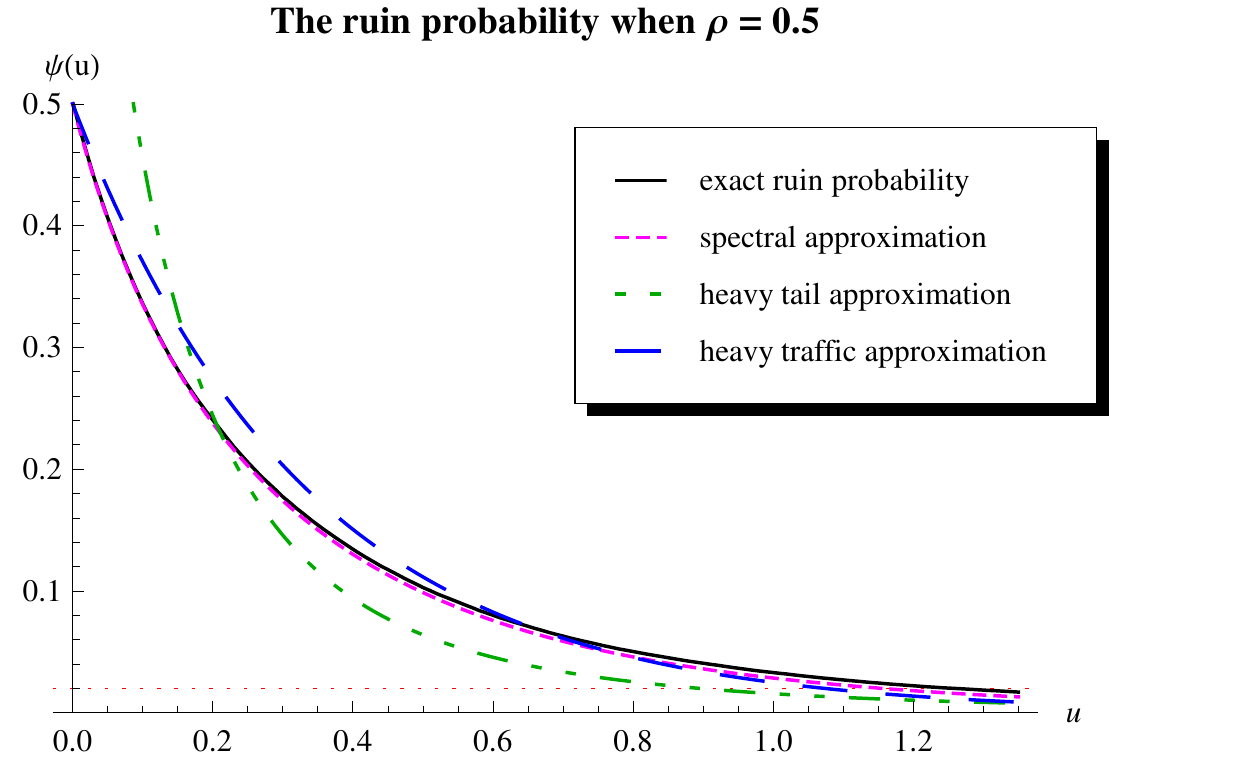}\label{figure:Pareto 0.5, with accuracy 0.02}
    }
    \subfigure[]
    {
        \includegraphics[scale=0.9]{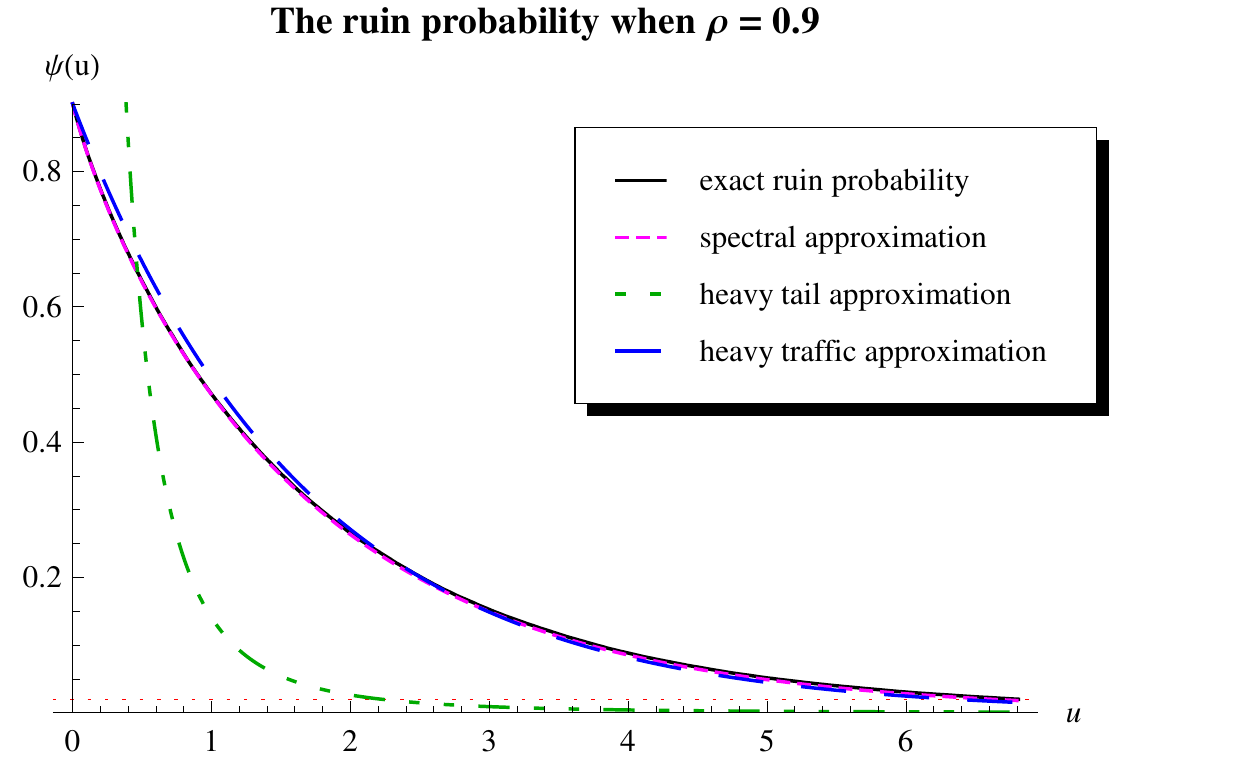}\label{figure:Pareto 0.9, with accuracy 0.02}
    }
\label{fig:subfigureExample}
\caption[]{The spectral approximation for guaranteed bound $\delta =0.02$, when the claims follow the Pareto$(4,3)$ distribution.}
\end{figure}

\begin{table}[ht!]
\vspace{5pt}
     \begin{center}
         \subtable[10 phases]
         {
             \begin{tabular}{|c|ccc|}\hline
                         $\rho$ & \text{bound}         & \text{max error}     & \text{ratio} \\ \hline
                         0.1    & 0.010                & 0.0048               & 2.11 \\
                         0.2    & 0.023                & 0.0106               & 2.13 \\
                         0.3    & 0.039                & 0.0180               & 2.17 \\
                         0.4    & 0.061                & 0.0275               & 2.21 \\
                         0.5    & 0.091                & 0.0401               & 2.27 \\
                         0.6    & 0.136                & 0.0580               & 2.35 \\
                         0.7    & 0.212                & 0.0849               & 2.50 \\
                         0.8    & 0.364                & 0.1299               & 2.80 \\
                         0.9    & 0.818                & 0.2263               & 3.61 \\\hline
              \end{tabular}
              \label{table:Abate ratios 10 phases}
         }
         \subtable[20 phases]
         {
             \begin{tabular}{|c|ccc|}\hline
                         $\rho$ & \text{bound}         & \text{max error}      & \text{ratio} \\ \hline
                         0.1    & 0.005                & 0.0026                & 2.06 \\
                         0.2    & 0.012                & 0.0057                & 2.08 \\
                         0.3    & 0.020                & 0.0097                & 2.09 \\
                         0.4    & 0.032                & 0.0150                & 2.12 \\
                         0.5    & 0.048                & 0.0222                & 2.15 \\
                         0.6    & 0.071                & 0.0326                & 2.19 \\
                         0.7    & 0.111                & 0.0490                & 2.27 \\
                         0.8    & 0.190                & 0.0787                & 2.42 \\
                         0.9    & 0.429                & 0.1479                & 2.90\\ \hline
              \end{tabular}
              \label{table:Abate ratios 20 phases}
           }
           \subtable[100 phases]
           {
              \begin{tabular}{|c|ccc|}\hline
                         $\rho$ & \text{bound}         & \text{max error}     & \text{ratio} \\ \hline
                         0.1    & 0.001                & 0.0005               & 2.02 \\
                         0.2    & 0.002                & 0.0012               & 2.02 \\
                         0.3    & 0.004                & 0.0021               & 2.02 \\
                         0.4    & 0.007                & 0.0033               & 2.03 \\
                         0.5    & 0.010                & 0.0049               & 2.04 \\
                         0.6    & 0.015                & 0.0073               & 2.05 \\
                         0.7    & 0.023                & 0.0112               & 2.06 \\
                         0.8    & 0.040                & 0.0189               & 2.10 \\
                         0.9    & 0.089                & 0.0406               & 2.19 \\ \hline
                \end{tabular}
                \label{table:Abate ratios 100 phases}
             }
     \end{center}
\caption{Ratios between the guaranteed bound and the maximum error of the spectral approximation, when the claims follow the Abate-Whitt distribution with $\mu=2$.}\label{table:Abate with ratios}
\end{table}

\begin{table}[ht!]
  \begin{center}
    \begin{tabular}{|c|ccccc|}\hline
      $\rho$  & \text{HT bound} & $k^*$    & \text{sp bound}       & \text{max HT error} & \text{max sp error} \\ \hline
        0.82  & 0.78               & 5     & 0.76                  & 0.0438              & 0.0312 \\
        0.85  & 0.65               & 8     & 0.63                  & 0.0403              & 0.0253 \\
        0.88  & 0.52               & 13    & 0.52                  & 0.0361              & 0.0196 \\
        0.91  & 0.39               & 25    & 0.39                  & 0.0304              & 0.0139 \\
        0.94  & 0.26               & 59    & 0.26                  & 0.0234              & 0.0081 \\
        0.97  & 0.13               & 248   & 0.13                  & 0.0144              & 0.0013 \\ \hline
    \end{tabular}
  \end{center}
  \caption{Comparison between the maximum heavy traffic and spectral errors, when the claims follow the Weibull$(0.5,3)$ distribution.}\label{table:Weibull comparison spectral with heavy traffic}
\end{table}

\begin{table}[ht!]
  \begin{center}
    \begin{tabular}{|c|ccccc|}\hline
         $\rho$ & \text{HT bound}    &$k^*$ & \text{sp bound}       & \text{max HT error} & \text{max sp error} \\ \hline
         0.82   & 0.90               & 4    & 0.91                  & 0.0392              & 0.0453 \\
         0.85   & 0.75               & 7    & 0.71                  & 0.0365              & 0.0387 \\
         0.88   & 0.60               & 11   & 0.61                  & 0.0326              & 0.0330 \\
         0.91   & 0.45               & 21   & 0.46                  & 0.0279              & 0.0261 \\
         0.94   & 0.30               & 51   & 0.30                  & 0.0226              & 0.0166 \\
         0.97   & 0.15               & 215  & 0.15                  & 0.0156              & 0.0074 \\ \hline
    \end{tabular}
  \end{center}
  \caption{Comparison between the maximum heavy traffic and spectral errors, when the claims follow the Pareto$(4,3)$ distribution.}\label{table:Pareto comparison spectral with heavy traffic}
\end{table}

\begin{table}[ptbh]
  \begin{center}
    \begin{tabular}{|c|ccccc|}\hline
         $\rho$ & \text{HT bound}   &$k^*$ & sp bound              & \text{max HT error} & \text{max sp error} \\ \hline
          0.82  & 0.568             & 7    & 0.569                 & 0.0051              & 0.0068 \\
          0.85  & 0.473             & 11   & 0.472                 & 0.0060              & 0.0066 \\
          0.88  & 0.379             & 18   & 0.386                 & 0.0044              & 0.0044 \\
          0.91  & 0.284             & 35   & 0.281                 & 0.0047              & 0.0026 \\
          0.94  & 0.190             & 82   & 0.189                 & 0.0047              & 0.0014 \\
          0.97  & 0.095             & 340  & 0.095                 & 0.0024              & 0.0025 \\ \hline
    \end{tabular}
  \end{center}
  \caption{Comparison between the maximum heavy traffic and spectral errors, when the claims follow the Pareto$(15.6,2.7)$ distribution.}\label{table:Pareto comparison spectral with heavy traffic 2}
\end{table}
\end{appendix}

\end{document}